
\input amstex

\magnification 1200
\loadmsbm
\parindent 0 cm

\define\nl{\bigskip\item{}}
\define\snl{\smallskip\item{}}
\define\inspr #1{\parindent=20pt\bigskip\bf\item{#1}}
\define\iinspr #1{\parindent=27pt\bigskip\bf\item{#1}}
\define\einspr{\parindent=0cm\bigskip}

\define\ot{\otimes}

\define\btr{\blacktriangleright}
\define\btl{\blacktriangleleft}
\define\dubbeldualA{\widehat{\widehat{\hbox{\hskip -0.13em $A$\hskip 0.13em }}}}
\define\dubbeldualepsilon{\widehat{\widehat{\hbox{\hskip -0.13em $\varepsilon$\hskip 0.13em }}}}
\define\dubbeldualphi{\widehat{\widehat{\hbox{\hskip -0.13em $\varphi$\hskip 0.13em }}}}

\input amssym
\input amssym.def

\centerline{\bf Algebraic Quantum Hypergroups}
\bigskip\bigskip
\centerline{\it L.\ Delvaux \rm ($^{*}$) and \it A.\ Van Daele \rm ($^{**}$)}
\bigskip\bigskip\bigskip
{\bf Abstract} 
\bigskip 
An algebraic quantum group is a regular multiplier Hopf algebra with integrals.  In this paper we will develop a theory of {\it algebraic quantum hypergroups}.  It is very similar to the theory of algebraic quantum groups, except that the comultiplication is no longer assumed to be a homomorphism.  We still require the existence of a left and of a right integral. There is also an antipode but it is characterized in terms of these integrals.  We construct the dual, just as in the case of algebraic quantum groups and we show that the dual of the dual is the original quantum hypergroup.  We define algebraic quantum hypergroups of compact type and discrete type and we show that these types are dual to each other. The algebraic quantum hypergroups of compact type are essentially the algebraic ingredients of the compact quantum hypergroups as introduced and studied (in an operator algebraic context) by Chapovsky and Vainerman. 
\snl
We will give some basic examples in order to illustrate different aspects of the theory. In a separate note, we will consider more special cases and more complicated examples. In particular, in that note, we will give a general construction procedure and show how known examples of these algebraic quantum hypergroups fit into this framework.
\nl
\nl
{\it April 2007} (Version 1.2)
\vskip 5 cm
\hrule
\bigskip
($^{*}$) Department of Mathematics, University of Hasselt, Agoralaan, B-3590 Diepenbeek, Belgium. E-mail: Lydia.Delvaux\@uhasselt.be
\snl
($^{**}$) Department of Mathematics, K.U.\ Leuven, Celestijnenlaan 200B,
B-3001 Heverlee (Belgium). E-mail: Alfons.VanDaele\@wis.kuleuven.be

\newpage

\bf 0. Introduction \rm
\nl
Let $A$ be a Hopf algebra (cf.\ [A] and [S]). The linear dual space $A'$ is made into an associative algebra if it is endowed with the product, dual to the coproduct on $A$. If $A$ is finite-dimensional, this dual algebra can be made into a Hopf algebra if the coproduct is defined, dual to the product in $A$. This nice duality breaks down when $A$ is no longer assumed to be finite-dimensional. In this case, the natural candidate for the coproduct on $A'$ will no longer map into the tensor product $A'\ot A'$, but rather in the (strictly bigger) algebra $(A\ot A)'$.
\snl
However, if we allow a more general structure, namely if we consider regular multiplier Hopf algebras with integrals (cf.\ [VD1] and [VD2]), we again have a nice duality. It generalizes the case of finite-dimensional Hopf algebras to a much greater class of objects. If e.g.\ $A$ is the group algebra of a group $G$, finite or not, then the dual exists in this more general framework and it is the multiplier Hopf algebra $K(G)$ of all functions on $G$ with finite support and with the coproduct, properly defined, dual to the product in $G$ (and the product in $A$). The class of these so-called {\it algebraic quantum groups}, contains the compact quantum groups, the discrete quantum groups, and many more cases.
\snl
In this paper, we will show that many of the nice aspects of the duality for these algebraic quantum groups remain valid if we consider algebras with comultiplications that are no longer assumed to be algebra homomorphisms. In this way, we naturally arrive at the study of {\it algebraic quantum hypergroups}. For a precise definition, we refer to Definition 1.10 in Section 1 of this paper. For algebraic quantum hypergroups, we get a duality, very much along the same lines as for the algebraic quantum groups (i.e.\ the regular multiplier Hopf algebras with integrals).
\snl
Hypergroups appear naturally when considering non-normal subgroups of a group. For simplicity, consider a {\it finite} group $G$ with a subgroup $H$. Denote by $A$ the algebra of complex functions on $G$ that are constant on double cosets. One can then define a linear map $\Delta: A \to A \ot A$ by
$$\Delta(f)(p,q)=\frac{1}{n}\sum_{h\in H} f(phq)$$
where $n$ is the number of elements of $H$ and $p,q$ are in $G$. This map will only be an algebra homomorphism if the subgroup $H$ is a normal subgroup. However, it will still satisfy coassociativity. A similar example can be constructed when $H$ is a finite subgroup of any group $G$, finite or not (see Example 1.11 in Section 1). We will use this example throughout the paper to motivate definitions and results.
\nl
It is well-known that the theory of algebraic quantum groups eventually led to a nice theory of locally compact quantum groups (cf.\ [K-V1], [K-V2] and [K-V3], see also [VD4] and references therein). Therefore, it is expected that the theory of algebraic quantum hypergroups, as developed in this paper, will serve as a source of inspiration for a possible theory of locally compact quantum hypergroups.
\snl
The compact quantum hypergroups, as developed by Chapovski and Vainerman in [C-V], and studied further (cf.\ e.g.\ [V] and [Ka]), should be a special case of these locally compact quantum hypergroups, just as the compact quantum groups, as developed by Woronowicz (cf.\ [W1] and [W2], see also [M-VD]) are a special case of the general locally compact quantum groups. In this paper we also define algebraic quantum hypergroups of compact type. These are essentially the compact quantum hypergroups of [C-V], but formulated in a purely algebraic context, just as the algebraic quantum groups of compact type (cf.\ [VD2]) are essentially the algebraic versions of the compact quantum groups of [W2]. We will not be able to prove this statement in a correct way. We plan to give a precise treatment in the future. However, in Section 5, we will shortly indicate why this result is to be expected and how it should be proven.
\snl
This is only one situation that supports the idea that the algebraic quantum hypergroups can be considered as an algebraic version of a forthcoming analytical theory of locally compact quantum hypergroups (including the compact quantum hypergroups of [C-V]). Another case to consider is found in the papers [L-VD1] and [L-VD2]. These papers deal with compact and discrete (quantum) subgroups of algebraic quantum groups. And whereas the first paper is of a purely algebraic nature, they both fit completely within the analytical theory. They are closely related. In fact, what motivated this paper on algebraic quantum hypergroups, are some results in [L-VD1]. We will come back to this relation when we draw conclusions in Section 5.
\nl
For these reasons, this paper should be of interest, not only for algebraists, but also for operator algebraists. The paper is written with this intention in mind. It is one of the reasons why we are also interested in the $^*$-algebra case and why sometimes we require positivity of the integrals.
\nl
{\it The paper is organized as follows}.
\nl
In {\it Section 1}, we give a precise definition of an algebraic quantum hypergroup. Essentially, it is like an algebraic quantum group, but without the requirement that the coproduct is an algebra homomorphism. We motivate this definition using a simple example (see earlier in this introduction). We show that the algebraic quantum groups, in the sense of [VD2], are special cases. We define the notion of integrals in this section. Also here, we define what we mean by an algebraic quantum group of compact type. 
\snl
In {\it Section 2} we prove that the integrals are unique, up to a scalar. We show that, just as in the case of algebraic quantum groups, left and right integrals are related by the so-called modular multiplier $\delta$ in $M(A)$.  In fact, many other features of the algebraic quantum groups are still valid, also in this more general setting. We obtain the same data and the same relations among these data.
\snl
In {\it Section 3}, we start with an algebraic quantum hypergroup $(A, \Delta)$ and we construct the dual algebraic quantum hypergroup $(\widehat{A}, \widehat{\Delta})$.  The main result is given in Theorem 3.11. In the next result, Theorem 3.12, we show that the dual of $(\widehat{A}, \widehat{\Delta})$ is again the original pair $(A,\Delta)$. We define here the algebraic quantum groups of discrete type and we show that these are dual to the ones of compact type, defined earlier in Section 1.
\snl
In {\it Section 4} we consider the different data, associated with the dual. They are given in terms of those of the original pair $(A, \Delta)$. We are able to obtain more relations between the different data for $(A,\Delta)$ and those of the dual $(\widehat{A}, \widehat{\Delta})$. We also look at the obvious module structures. Moreover, we recover e.g.\ a formula for the fourth power of the antipode, which is known in the theory of Hopf algebras as Radford's formula. We see again that many of the results, known for algebraic quantum groups stay valid for the algebraic quantum hypergroups.
\snl
In {\it Section 5}, we conclude and we discuss some open problems and give some ideas for further research. One of the obvious things to do is to see whether, as in the case of algebraic quantum groups, it is possible also here to start with a $^*$-algebraic quantum hypergroup with positive integrals and lift it to what should be called a C$^*$-algebraic quantum hypergroup (as in [K-VD]). Then it should be possible to give a precise proof of the statement that the compact quantum hypergroups, introduced in [C-V] are in correspondence with $^*$-algebraic quantum hypergroups with positive integrals and of compact type as they are defined in this paper.
\snl
In a forthcoming paper [D-VD2] we plan to look at more (and more complicated) examples, as well as the study of some construction methods. 
\nl
In this paper, we will work with the following {\it assumptions} and {\it notations}. All algebras are associative and are considered over the field $\Bbb C$. We do not assume that an algebra $A$ is unital, but we require that the multiplication on $A$, considered as a bilinear map, is non-degenerate.
\snl
We can consider the multiplier algebra $M(A)$ of an algebra $A$, see e.g.\ the appendix in [VD1].  It can be characterized as the largest unital algebra which contains $A$ as a dense two-sided ideal. It is easy to see that algebra (anti)-isomorphisms extend in a unique way to the multiplier algebras.  More generally, for two algebras $A$ and $B$, let $\alpha : A \to M(B)$ be an algebra homomorphism. Then $\alpha$ is called {\it non-degenerate} if $\alpha (A) B = B \alpha(A) = B$.  It is possible to extend $\alpha$, in a natural way, to a unital homomorphism from $M(A)$ to $M(B)$.  This extension is unique and is again denoted by the symbol $\alpha$, see [VD1, Proposition A5].  The identity in various algebras is denoted by 1.  The symbol $\iota$ denotes the identity map.  The linear dual of an algebra $A$ is denoted as $A'$.  For $\omega \in A'$, we consider the slice maps $\iota \otimes \omega: A \otimes A \to A$ and $\omega \otimes \iota: A \otimes A \to A$.  For a subset $W \subset A$ we denote by $\text{sp}(W)$ the linear space generated by the elements of $W$.
\nl
\bf Acknowlegdements \rm
\snl
The second author would like to thank M.B.\ Landstad of the University of Trondheim (Norway) for his hospitaly during the many visits and for the fruitful cooperation. This paper grew out of the joint work on discrete and compact subgroups of algebraic quantum groups ([L-VD1]).  

\newpage

\bf 1. The definition of an algebraic quantum hypergroup \rm
\nl
In this section, we will develop the notion of an algebraic quantum hypergroup (see Definition 1.10). We will illustrate various aspects of this definition by means of a simple, but typical example (refered to as the motivating example further in the paper, see Example 1.11). We will also explain that an algebraic quantum hypergroup, as defined in this section, is very much like an algebraic quantum group, but without assuming that the coproduct is an algebra homomorphism (see Proposition 1.14 in this section and Proposition 2.3 in the next one). The main properties of an algebraic quantum hypergroup will be obtained in the next section. We will see that many of the properties of an algebraic quantum group remain true for these quantum hypergroups. This is somewhat remarkable.
\nl 
Our starting point is an (associative) algebra $A$ over $\Bbb C$, with or without identity, but with a non-degenerate product.  The tensor product $A \otimes A$ of $A$ with itself is  again an algebra in the obvious way and the product is still non-degenerate.  We can consider the multiplier algebras $M(A)$ and $M(A \otimes A)$ of $A$ and $A \otimes A$ respectively.  As usual, we view $M(A) \otimes M(A)$ as sitting in $M(A \otimes A)$.
\snl
We start with the following definition.

\inspr{1.1} Definition \rm 
Let $A$ be as above.  A {\it comultiplication} (or coproduct) on $A$ is a linear map $\Delta : A \rightarrow M(A \otimes A)$ such that both $\Delta(a) (1 \otimes b)$ and $(a \otimes 1) \Delta(b)$
belong to $A \otimes A$ for all $a, b \in A$ and such that 
$$(a \otimes 1 \otimes 1) (\Delta \otimes \iota)(\Delta(b) (1 \otimes c)) = (\iota \otimes \Delta)((a \otimes 1) \Delta(b)) (1 \otimes 1 \otimes c)$$
for all $a, b, c \in A$.
\einspr

Recall that we use $1$ for the identity in $M(A)$ and $\iota$ for the identity map on $A$.
\snl
If $A$ has an identity, then $M(A) = A$ and $M(A \otimes A) = A \otimes A$ and then a linear map $\Delta : A \to A \otimes A$ is a comultiplication when it is coassociative.  Therefore, in general, the last condition in the definition will also be called coassociativity of $\Delta$.  The first condition is needed in order to be able to formulate the second one.
\snl
Observe that in this paper we do not assume that $\Delta$ is an algebra homomorphism.
\snl
If $A$ is a $^\ast$-algebra, we will assume also that $\Delta$ is a $^\ast$-map, i.e.\ we require that $\Delta(a^\ast) = \Delta(a)^\ast$ where the involution on $M(A \otimes A)$ comes in a natural way from the involution on $A \otimes A$.

\inspr{1.2} Definition \rm 
Let $\Delta$ be a comultiplication on $A$ as in Definition 1.1.  Then $\Delta$ is called {\it regular} if also 
$\Delta(a) (b \otimes 1)$ and $(1 \otimes a) \Delta(b)$ are in $A \otimes A$ for all $a, b \in A$.
\einspr

In the case of a $^\ast$-algebra, regularity of the coproduct is automatic.
\snl
It is not hard to show that for a regular comultiplication $\Delta$, the opposite map $\Delta'$, obtained from $\Delta$ by composing it with the flip $\zeta:A \otimes A \rightarrow A \otimes A$ (defined by $\zeta (a \otimes b) = b \otimes a$ and extended to $M(A \otimes A)$, will again be coassociative and hence a comultiplication.  
\snl
Further in this paper, we will only work with regular comultiplications.
\nl
An algebraic quantum hypergroup will be defined as a pair $(A,\Delta)$ where $A$ is an algebra (with a non-degenerate product) and $\Delta$ a regular comultiplication satisfying certain extra assumptions (see Definition 1.10 below). We develop the precise definition in a few steps. 
\snl
First we assume the {\it existence of a counit} as in the following definition.

\inspr{1.3} Definition \rm 
Let $(A, \Delta)$ be an algebra $A$ with a coproduct $\Delta$.  A homomorphism $\varepsilon: A \rightarrow \Bbb C$ is called a {\it counit} if
$(\varepsilon \otimes \iota) \Delta(a) = a$ and $(\iota \otimes \varepsilon) \Delta (a) = a$ for all $a \in A$.
\einspr

The above definition makes sense because, for a regular comultiplication, $(\omega \otimes \iota)\Delta(a)$ and $(\iota \otimes \omega)\Delta(a)$ can be defined in $M(A)$ for all $\omega \in A'$ and $a \in A$.  The first condition in the definition means e.g.\ that for all $a, b \in A$ we have
$$(\varepsilon \otimes \iota) (\Delta(a)(1\otimes b)) = ab$$
(and similarly for the 3 other possibilities).
\snl
Now, we show that the counit, if it exists, must be unique.  In fact, we get a slightly stronger result.

\inspr{1.4} Proposition \rm 
Let $(A, \Delta)$ be as before and assume that $\varepsilon$ is a counit.  If $\varepsilon'$ is any linear map from $A$ to $\Bbb C$ satisfying $(\iota \otimes \varepsilon')\Delta(a) = a$ for all $a \in A$, then $\varepsilon' = \varepsilon$.  Similarly, if $(\varepsilon' \otimes \iota) \Delta(a) = a$ for all $a \in A$, then we also have $\varepsilon' = \varepsilon$.
\snl\bf Proof: \rm
Assume that $\varepsilon' \in A'$ and that $(\iota \otimes \varepsilon')\Delta(a) = a$ for all $a \in A$.  We know that $(\varepsilon \otimes \iota) ((b\otimes 1)\Delta(a)) = \varepsilon(b) a$ because $\varepsilon$ is a counit (and so a homomorphism).  If we apply $\varepsilon'$ to this equation we get $\varepsilon(ba) = \varepsilon(b) \varepsilon'(a)$. Because $\varepsilon(ba) = \varepsilon(b) \varepsilon(a)$ and $\varepsilon$ can not be trivially zero, we get $\varepsilon(a) = \varepsilon'(a)$ for all $a \in A$.  This proves the first statement.  The other one is obtained in a similar way (or by symmetry). \hfill $\blacksquare$
\einspr

If $A$ is a $^\ast$-algebra, we will have that $\varepsilon$ is a $^\ast$-homomorphism.  Indeed, one can easily verify that $\varepsilon_1 : A \rightarrow \Bbb C$ defined by $\varepsilon_1(a) = \varepsilon(a^\ast)^-$ (where $\overline{\lambda}$ is the complex conjugate of $\lambda \in \Bbb C$), is again a counit and so $\varepsilon_1 = \varepsilon$.  This means that $\varepsilon$ is a $^\ast$-homomorphism.
\snl
From now on, we assume that $(A, \Delta)$ is a pair of an algebra $A$ with a regular coproduct $\Delta$ and we assume that a counit $\varepsilon$ exists.
\nl
Next, it would be most natural to introduce the notion of an antipode.  However, in this setting, it turns out to be appropriate to first consider a (left) integral. The reader may compare this with the Larson-Sweedler theorem for multiplier Hopf algebras (cf.\ [VD-W]) where the antipode is proven from the existence of the integrals.
\snl
As expected, we have the following definition for integrals in this setting.

\inspr{1.5} Definition \rm 
A non-zero linear functional $\varphi$ on $A$ is called a {\it left integral} if
$$(\iota \otimes \varphi) \Delta(a) = \varphi(a) 1$$
in $M(A)$ for all $a \in A$. Similarly, a non-zero linear functional $\psi$ is called a right integral if
$$(\psi \otimes \iota) \Delta(a) = \psi (a) 1$$
in $M(A)$ for all $a \in A$.
\einspr

If $A$ is a $^\ast$-algebra we will (at least) assume that $\varphi$ is self-adjoint, i.e. that $\varphi(a^\ast) = \varphi(a)^-$ for all $a \in A$.  This in fact is not really an extra assumption. Indeed, if a left integral $\varphi$ exists, then there must also exist a self-adjoint one. Take $\varphi + \overline{\varphi}$ or $i(\varphi-\overline{\varphi})$ where $\overline{\varphi}$ is the left integral on $A$ defined as $\overline{\varphi}(a) = \varphi(a^\ast)^-$ for all $a \in A$.  It makes sense to assume that there is a positive left integral (i.e. $\varphi(a^\ast a) \geq 0$ for all $a \in A$).  This however is not a trivial assumption.
\nl
With the assumptions we have made so far, that is having a regular comultiplication and a left integral, we can prove already, just as in the case of algebraic quantum groups (see Proposition 2.6 in [Dr-VD-Z]), that the algebra $A$ must have {\it local units} (in the sense of the following proposition). The argument is completely the same as in the case of algebraic quantum groups, i.e.\ there is no need for the coproduct to be a homomorphism.

\inspr{1.6} Proposition \rm 
Let $(A, \Delta)$ be as before (in particular, we assume the existence of a left integral $\varphi$). Given elements $\{a_1, a_2,\ldots, a_n\}$, there exists an element $e \in A$ such that $a_ie = ea_i = a_i$ for all $i$.
\snl \bf Proof: \rm
Define the linear space $V$ in $A^{2n}$ as follows
$$V = \{(aa_1, aa_2,\ldots, aa_n, a_1a, a_2a,\ldots,a_na)\mid a \in A\}.$$
Consider a linear functional on $A^{2n}$ that is zero on $V$.  This means that we have functionals $\omega_i$ and $\rho_i$ on $A$ for $i=1,\ldots,n$, such that
$$\sum\limits_{i=1}^n \omega_i (aa_i) + \sum\limits_{i=1}^n \rho_i(a_ia) = 0$$
for all $a \in A$. Then, for all $x, a \in A$ we have
$$\align
x &\left(\sum\limits_{i=1}^n (\omega_i \otimes \iota)(\Delta(a) (a_i \otimes 1)) + \sum\limits_{i=1}^n (\rho_i \otimes \iota) ((a_i \otimes 1)\Delta(a))\right)\\
&\qquad = \sum\limits_{i=1}^n (\omega_i \otimes \iota)((1 \otimes x)\Delta(a) (a_i \otimes 1)) + \sum\limits_{i=1}^n (\rho_i \otimes \iota) ((a_i \otimes x) \Delta(a)) = 0.
\endalign$$
As the product in $A$ is non-degenerate, we get for all $a \in A$ that
$$\sum\limits_{i=1}^n (\omega_i \otimes \iota)(\Delta(a) (a_i \otimes 1)) + \sum\limits_{i=1}^n (\rho_i \otimes \iota) ((a_i \otimes 1) \Delta(a)) = 0.$$
If we now apply $\varphi$ on this expression, we obtain
$$\varphi(a) \left(\sum\limits_{i=1}^n \omega_i (a_i) + \sum\limits_{i=1}^n \rho_i (a_i)\right) = 0.$$
for all $a\in A$. Therefore, as $\varphi$ is non-zero, we may conclude $\sum\limits_{i=1}^n \omega_i (a_i) +  \sum\limits_{i=1}^n \rho_i(a_i) = 0$.
\snl
So, any linear functional on $A^{2n}$ that is zero on the space $V$ is also zero on the vector $(a_1, a_2,\ldots, a_n, a_1, a_2,\ldots,a_n)$. Therefore, $(a_1, a_2, \ldots, a_n, a_1, a_2,\ldots, a_n)$ belongs to the space $V$.   This means that there exists an element $e \in A$ such that $ea_i = a_i$ and $a_ie = a_i$ for all $i$.\hfill $\blacksquare$
\einspr

This result has an important practical consequence. We formulate it as a separate remark. 
 
\inspr{1.7} Remark \rm  
For regular multiplier Hopf algebras, the use of the Sweedler notation is justified, just as for ordinary Hopf algebras, see e.g.\ [Dr-VD-Z].  Also for a pair $(A, \Delta)$ with an algebra $A$ having local units as in Proposition 1.6, we can use a formal expression for $\Delta(a)$ when $a \in A$. Even though $\Delta(a)$ for $a\in A$ in general is not in $A \otimes A$, we do have that $\Delta(a) (1 \otimes b) \in A  \otimes A$ for all $b \in A$. By Proposition 1.6, we know that there is an element $e \in A$ such that $eb = b$.  Therefore, $\Delta(a) (1 \otimes b) = \Delta(a) (1 \otimes e) (1 \otimes b)$ and we can write $\Delta(a) (1 \otimes b) = \sum a_{(1)} \otimes a_{(2)} b$.  The expression $\sum a_{(1)} \otimes a_{(2)}$ stands for $\Delta(a) (1 \otimes e)$.  Observe that this expression is dependent on the element $b$, but for several elements $b_i$ we can choose the same element $e$ (and so we can use the same expression).
\snl
The {\it Sweedler notation} in this context has to be used with some care, but it is very convenient to make formulas more transparant. Further in this paper, we will indeed, when appropriate, make use of the Sweedler notation (in the above sense). Observe however that it is always possible to translate formulas in such a way that the Sweedler notation is not needed. 
\einspr

\snl
Without further assumptions, not much can be proven about integrals in this context.  We will need to impose an extra condition on $\varphi$.  One condition is faithfulness. Recall that a linear functional $f$ on $A$ is called {\it faithful} if, for $a \in A$, we must have $a = 0$ when either $f(ab) = 0$ for all $b \in A$ or $f(ba) = 0$ for all $b \in A$.
\snl
We will also need a condition expressing the existence of an antipode.  Before we can do this, we require the following lemma.

\inspr{1.8} Lemma \rm 
Let $(A,\Delta)$ be as before.  In particular, we assume the existence of a counit $\varepsilon$.  If $f$ is a faithful linear functional on $A$, then for any $a \in A$, there is an element $e \in A$ such that
$$a = (\iota \otimes f) (\Delta (a) (1 \otimes e)).$$
\snl \bf Proof: \rm
Take $a \in A$ and define $V = \{(\iota \otimes f)(\Delta (a) (1 \otimes b)) \mid b \in A\}$. We need to show that $a \in V$.  Suppose that this is not the case.  Then, there is an element $\omega \in A'$ such that $\omega(a) \neq 0$ while $\omega|_V = 0$.  This last property means however that $f(xb)= 0$ for all $b \in A$ where $x = (\omega \otimes \iota) \Delta(a)$. Observe that $x \in M(A)$ and not necessarily $x \in A$.  However, we get $f(xb'b'') = 0$ for all $b', b'' \in B$ and by the faithfulness of $f$, we must have  $xb' = 0$ for all $b' \in A$.  If we apply the counit, we get $\omega(a) \varepsilon(b') = 0$ and hence $\omega(a) = 0$.   This is a contradiction. \hfill $\blacksquare$
\einspr

In a completely similar way we get
$$\align
&a \in \{(\iota \otimes f) ((1\otimes b) \Delta(a)) \mid b \in A\}\\
&a \in \{(f \otimes \iota)((b \otimes 1) \Delta(a)) \mid b \in A\}\\
&a \in \{(f \otimes \iota)(\Delta(a) (b \otimes 1)) \mid b \in A\}
\endalign$$
for any faithful $f \in A'$. In particular, when we assume that a left integral $\varphi$ is faithful, we have
$$\align
&A = \text{sp} \{(\iota \otimes \varphi)(\Delta (a) (1 \otimes b)) \mid a, b \in A\}\\
&A = \text{sp} \{(\iota \otimes \varphi)((1\otimes a)\Delta(b)) \mid a, b \in A\}
\endalign$$
where $\text{sp}$ is used to denote the linear span of a set of elements in $A$.
\einspr

Now we are ready to consider the existence of an antipode.

\inspr{1.9} Definition \rm 
Let $(A,\Delta)$ be as before and assume that there is a faithful left integral $\varphi$.  Suppose that there is a linear bijective map $S : A \rightarrow A$ satisfying 
$$S((\iota \otimes \varphi)(\Delta(a)(1\otimes b))) = (\iota \otimes \varphi)((1 \otimes a) \Delta(b))$$
for all $a, b \in A$.  Observe that, as a consequence of the previous lemma, this linear map, when it exists, is uniquely determined by the above formula. If moreover this map $S$ is a anti-homomorphism, then $S$ is called the {\it antipode} (relative to $\varphi)$.
\einspr

If $A$ is a $^\ast$-algebra and if $\varphi$ is self-adjoint, then we get $S(x)^\ast = (\iota \otimes \varphi) (\Delta(b^\ast)(1\otimes a^\ast))$ when $x = (\iota \otimes \varphi) (\Delta(a)(1\otimes b))$ and we see that $S(S(x)^\ast)^\ast = x$ for all $x \in A$.
\nl
Later, we will show that left integrals are unique (provided there exists a left integral $\varphi$ with an antipode relative to this integral $\varphi$); see Proposition 2.4 in the next section.  This is why we can now formulate the following main definition.

\iinspr{1.10} Definition \rm 
Let $(A,\Delta)$ be an algebra with a regular comultiplication $\Delta$ and a counit.  Assume that there exists a faithful left integral $\varphi$ with  an antipode $S$ (relative to $\varphi$).  Then $(A,\Delta)$ is called an {\it algebraic quantum hypergroup}.  If moreover $A$ is a $^\ast$-algebra and $\Delta$ is a $^\ast$-map, then we call $(A,\Delta)$ a {\it $^\ast$-algebraic quantum hypergroup}.
\einspr

Before we continu with proving the first elementary properties, let us consider the following, motivating example.

\iinspr{1.11} Example \rm
Let $G$ be a (discrete) group and let $H$ be a finite subgroup. Let $A$ be the space of complex functions on $G$, with finite support and constant on double cosets of $G$ w.r.t.\ the subgroup $H$. So, $f(hph')=f(p)$ for all $h,h'\in H$ and $p\in G$ when $f\in A$. We need $H$ finite because we want $f$ to have finite support. It is clear that $A$ is an algebra when the product is defined pointwise. This product is non-degenerate. It becomes a $^*$-algebra if we define $f^*(p)=f(p)^-$ (where as before, $\lambda^-$ denotes the complex conjugate of a number $\lambda\in\Bbb C$).
\snl
Now define $\Delta$ by
$$\Delta(f)(p,q)=\frac1n \sum_{h\in H} f(phq)$$
where $n$ is the number of elements in $H$ and $p,q\in G$ and $f\in A$. It is not hard to verify that $\Delta(f)(1\otimes g)$ and $\Delta(f)(g\otimes 1)$ belong to $A\otimes A$. Also coassociativity is satisfied. We get e.g.\
$$((\Delta\otimes\iota)\Delta(f))(p,q,r)=\frac{1}{n^2}\sum_{h,h'\in H}f(phqh'r)$$
when $p,q,r\in G$ and $f\in A$.
So, $\Delta$ is a regular comultiplication on $A$ in the sense of Definitions 1.1 and 1.2 (observe that $A$ is abelian).
\snl
If we set $\varepsilon(f)=f(e)$ for $f\in A$ where $e$ denotes the identity element in $G$, we clearly get a counit (in the sense of Definition 1.3).
\snl
If we put $\varphi(f)=\sum_{p\in G} f(p)$, which is possible for $f\in A$ (as we have functions with finite support)), then we find
$$\frac1n\sum_{q\in G, h\in H} f(phq)= \sum_{q\in G}f(pq)=\sum_{q\in G} f(q)$$
and we see that $(\iota\otimes\varphi)\Delta(f)=\varphi(f)1$ for all $f\in A$. Therefore, we get a left integral in the sense of Defnition 1.5. In this case, $\varphi$ is also a right integral. It is faithful and positive (when the $^*$-algebra structure is considered).  
\snl
Finally, if $S$ is defined by $S(f)(p)=f(p^{-1})$ when $p\in G$ and $f\in A$, we see that indeed, $S$ is an antipode (relative to $\varphi$). If e.g.\ we have elements $f,g\in A$ and $p\in G$, we find
$$\align ((\iota\otimes\varphi)(\Delta(f)(1\otimes g)))(p)
       &= \frac1n \sum_{q\in G,\, h\in H} f(phq)g(q) \\
       &= \frac1n \sum_{q\in G,\, h\in H} f(pq)g(h^{-1}q) \\
       &= \sum_{q\in G} f(pq)g(q)
\endalign$$
(where we have used that $g$ is constant on cosets). Similarly, 
$$\align ((\iota\otimes\varphi)(1\otimes f)\Delta(g))(p)
       &= \sum_{q\in G} f(q)g(pq) \\
       &= \sum_{q\in G} f(p^{-1}q)g(q)
\endalign$$
and so $S$ indeed satisfies the required formula.
\snl
If we combine all the previous results, we see that $(A,\Delta)$ is an algebraic quantum hypergroup in the sense of Definition 1.10.
\einspr

Let us make a few more remarks w.r.t.\ this example.
\snl
If the subgroup $H$ is trivial, and only consists of the identity, then the above algebraic quantum hypergroup is simply the algebraic quantum group $K(G)$ of all complex functions with finite support on $G$ and the natural comultiplication. Observe that in this case, the coproduct is nothing else but the usual coproduct on $K(G)$. If (more generally) the subgroup $H$ is a normal subgroup, then we get the algebraic quantum group $K(G/H)$. Indeed, in this case
$$\Delta(f)(p,q)=\frac1n\sum_{h\in H}f(phq)=\frac1n\sum_{h\in H}f(pq(q^{-1}hq))=f(pq)$$
when $f\in A$ and $p,q\in G$. In fact, only if the group $H$ is normal, the coproduct is a homomorphism (and so only in this case we have an algebraic quantum group and not just an algebraic quantum hypergroup).
\snl
It is also instructive to illustrate Proposition 1.6 and Lemma 1.8 for this example. When $f$ is a function in $A$, it has a finite support which is a finite union of double $H$-cosets. If now $g$ is the function that is $1$ on this support and $0$ everywhere else, we find $g\in A$ and $fg=f$. Similarly, when we have a finite number of functions, $f_1, f_2, \ldots, f_m$, we simply consider the union of these supports. This illustrates 1.6. To illustrate 1.8, consider the function $g$ which is $\frac1n$ on $H$ (with $n$ the number of elements in $H$) and $0$ everywhere else. Again $g\in A$ and now, when $f\in A$ and $p,q\in G$ we have
$$(\Delta(f)(1\otimes g))(p,q)=\frac1n\sum_{h\in H}f(phq)g(q)=f(p)g(q)$$
so that $\Delta(f)(1\otimes g)=f\otimes g$. When we apply the left integral $\varphi$, we get $(\iota\otimes \varphi)(\Delta(f)(1\otimes g))=f$. We see that we can take the same $g$ for all $f$ in this case. This is so because we have a left co-integral (cf.\ Section 3). 
\nl
Apart from this example, we also look briefly at a special case. First, let us prove the following simple result.

\iinspr{1.12} Lemma \rm
Let $(A,\Delta)$ be an algebraic quantum hypergroup. If $A$ has an identity, then $\Delta(1)=1\otimes 1$.
\snl\bf Proof: \rm
Consider the antipode property with $b=1$ and $a$ arbitrary. We get
$$(\iota\ot\varphi)((1\otimes a)\Delta(1))=S((\iota\otimes\varphi)\Delta(a))=\varphi(a)S(1)=\varphi(a)1.$$
Because this holds for all $a$ and as $\varphi$ is faithful, we must have $\Delta(1)=1\otimes 1$.
\einspr

This leads to the following definition.

\iinspr{1.13} Definition \rm
An algebraic quantum hypergroup $(A,\Delta)$ is called of {\it compact type} if the algebra $A$ has an identity (and hence that $\Delta(1)=1\otimes 1$).
\einspr

The basic example given in 1.11 is of compact type, only if $G$ is finite. In general, it is of discrete type (cf.\ Definition 3.14 in Section 3). However, in Section 2 of [L-VD1], we encounter natural examples of algebraic quantum hypergroups of compact type. They are constructed from a so-called group-like projection in an algebraic quantum group. In the case of a group $G$, with $A$ being the algebra of complex functions on $G$ with finite support, considered with the pointwise product, such a group-like projection is obtained when we have a finite subgroup $H$ of $G$ and when we take the function that is $1$ on $H$ and $0$ anywhere else. Unfortunately, in this case, the resulting algebraic quantum group of compact type is noting else but the finite-dimensional Hopf algebra of all complex functions on $H$ with pointwise product and coproduct dual to the product in $H$. Only in the non-abelian case, non-trivial examples of compact type algebraic quantum hypergroups are found in this way. Again see Section 2 in [L-VD1] and also [D-VD2].
\snl
In Section 5, we will discuss more about these special cases and examples and also explain the terminology.

\nl
Let us now finish this section with the following important result.

\iinspr{1.14} Proposition \rm
If $(A,\Delta)$ is an algebraic quantum group (in the sense of [VD2]), then it is also an algebraic quantum hypergroup (in the sense of Definition 1.10 above).
\snl\bf Proof: \rm 
By assumption, $A$ is an algebra with a non-degenerate product and $\Delta$ is a (regular) coproduct in the sense of Definition 1.1 (it is even an algebra homomorphism in this case). For algebraic quantum groups, we have the existence of a counit (in the sense of Definition 1.3). We also have a faithful left integral $\varphi$ as in 1.5. There is an antipode $S$. It is a bijective, anti-isomorphism. The formula 
$$S((\iota \otimes \varphi)(\Delta(a)(1\otimes b))) = (\iota \otimes \varphi)((1 \otimes a) \Delta(b)),$$
needed to satisfy the requirements for an antipode in Definition 1.9, for all $a,b\in A$, is found in the proof of Proposition 3.11 in [VD2]. Therefore, all the assumptions in Definition 1.10 are satisfied and we do have an algebraic quantum hypergroup.
\einspr

Conversely, it is also true that an algebraic quantum hypergroup $(A,\Delta)$ with a coproduct $\Delta$ that is an algebra homomorphism, is actually an algebraic quantum group. Before we can show this however, we first need some basic properties of the antipode and we have chosen to prove these in the beginning of the next section.
\nl\nl

\bf 2. First properties of algebraic quantum hypergroups \rm
\nl
In this section, we consider an algebraic quantum hypergroup $(A,\Delta)$ as in Defintion 1.10 in the previous section. We will prove various properties, very similar as in the case of ordinary algebraic quantum groups. In particular, we will prove uniqueness of the integrals, we will get the scaling constant $\tau$, we will obtain the modular element $\delta$ relating the left with the right integral, we will get the modular automorphisms $\sigma$ and $\sigma'$ and we will obtain the various formulas relating these objects. The proofs are not more difficult than in the case of algebraic quantum groups. In fact, the way it is done here, although similar as in [VD2], is somewhat simpler.
\snl   
So consider in what follows an algebraic quantum hypergroup $(A,\Delta)$  with counit $\varepsilon$ and a faithful left integral $\varphi$ such that there exists an antipode $S$ relative to $\varphi$.
\snl
First we need some basic properties of the antipode. 

\inspr{2.1} Proposition \rm 
We have $\varepsilon (S(x)) = \varepsilon(x)$ for all $x\in A$. Also $\Delta (S(x)) = \zeta(S \otimes S) \Delta(x)$ whenever $x\in A$ where $\zeta$ is the flip on $A\otimes A$ (extended to $M(A \otimes A)$). We also use the extension of $S \otimes S$ to $M(A \otimes A)$ which is possible because $S$ is assumed to be a anti-isomorphism.

\snl\bf Proof: \rm 
Take $a, b \in A$ and set $x =(\iota \otimes \varphi)(\Delta(a)(1 \otimes b))$. By the definition of $S$, we have $S(x) = (\iota \otimes \varphi)((1 \otimes a)\Delta(b))$. If we apply $\varepsilon$ to both equations, we get $\varepsilon(x) = \varphi(ab)$ and $\varepsilon(S(x)) = \varphi(ab)$.  Because all elements in $A$ are of the form above, we have proven the first statement.
\snl
To prove the second statement, let $c$, $d \in A$. Then we have
$$\align
(c \otimes d) \Delta(S(x)) &= (\iota \otimes \iota \otimes \varphi) ((c \otimes d \otimes 1) (\Delta \otimes \iota)((1 \otimes a) \Delta(b)))\\
&= (\iota \otimes \iota \otimes \varphi)((1 \otimes d \otimes a)(\iota \otimes \Delta) ((c \otimes 1) \Delta(b)))\\
&= (1 \otimes d) (\iota \otimes S)(\iota \otimes \iota \otimes \varphi) ((c \otimes 1 \otimes 1) \Delta_{23}(a) \Delta_{13}(b))\\
&= (c \otimes 1) \zeta (S \otimes \iota)(\iota \otimes \iota \otimes \varphi)(\Delta_{13}(a) \Delta_{23}(b) (S^{-1}(d) \otimes 1 \otimes 1)).
\endalign$$
Observe that in the formulas above, we use the leg-numbering in the usual sense (e.g.\ $\Delta_{23}(a)$ means that $\Delta(a)$ is seen as multiplying  the second and the third components and leaving the first component fixed). If we cancel $c$ in the above formula, we can continu and for all $d \in A$ we have
$$\align
(1 \otimes d) \Delta(S(x)) &= \zeta (S \otimes \iota)(\iota \otimes \iota \otimes \varphi) (\Delta_{13}(a) \Delta_{23}(b) (S^{-1}(d) \otimes 1 \otimes 1))\\
&= \zeta (S \otimes S)(\iota \otimes  \iota \otimes \varphi) ((\Delta\otimes \iota) (\Delta(a) (1 \otimes b)) (S^{-1}(d) \otimes 1 \otimes 1))\\
&= \zeta(S \otimes S) (\Delta(x) (S^{-1}(d) \otimes 1))\\
&= \zeta((d \otimes 1) (S \otimes S) \Delta(x))\\
&= (1 \otimes d) \zeta (S \otimes S) (\Delta(x)).
\endalign$$
This means $\Delta(S(x)) = \zeta(S \otimes S) \Delta(x)$ in $M(A  \otimes A)$.  This completes the proof.
\hfill $\blacksquare$
\einspr

Also the next result is not unexpected.

\inspr{2.2} Proposition \rm 
Define $\psi = \varphi \circ S$.  Then $\psi$ is a faithful right integral on $A$.  Moreover
$$S((\psi \otimes \iota)((b \otimes 1) \Delta(a))) = (\psi \otimes \iota) (\Delta(b) (a \otimes 1))$$
for all $a,b\in A$.
\snl \bf Proof: \rm 
Because $S$ is a bijective anti-isomorphism and $\varphi$ is faithful, it is quite obvious that $\psi$ will again be faithful.  The right invariance of $\psi$ follows from the left invariance of $\varphi$ and because $S$ flips the coproduct.  To prove the formula relating $S$ and $\psi$, take $a, b \in A$ and start with the following equation 
$$(1)\qquad\qquad(\iota \otimes \varphi)(\Delta(a)(1 \otimes b)) = S^{-1} (\iota \otimes \varphi)((1 \otimes a) \Delta(b))$$
Using again that $S$ flips the coproduct (as proven in Proposition 2.1), we can write the left hand side of the equation (1) above as 
$$\align 
S(\iota \otimes \varphi)(S^{-1}\otimes \iota) (\Delta(a)(1 \otimes b)) &= S(\iota \otimes \psi)(S^{-1} \otimes S^{-1}) (\Delta(a) (1 \otimes b))\\
&= S(\psi \otimes \iota)((S^{-1} (b) \otimes 1) \Delta(S^{-1}(a))).
\endalign$$
The right hand side of the  equation (1) can be written as
$$(\iota \otimes \psi)(S^{-1} \otimes S^{-1}) ((1 \otimes a)\Delta(b)) = (\psi \otimes \iota) (\Delta(S^{-1}(b)) (S^{-1} (a) \otimes 1)).$$
Therefore, replacing $S^{-1}(a)$ by $c$ and $S^{-1}(b)$ by $d$, the equation (1) yields
$$S(\psi \otimes \iota)((d \otimes 1)\Delta(c)) = (\psi \otimes \iota)(\Delta(d) (c\otimes 1))$$
for all $c, d\in A$. This proves the proposition.\hfill $\blacksquare$
\einspr

Observe that $\varphi' \circ S$ (and of course also $\varphi' \circ S^{-1}$) will be a right integral for any left integral $\varphi'$.
\snl
Now, we are ready to prove another important result, announced already at the end of the previous section. 

\inspr{2.3} Proposition \rm 
Let $(A,\Delta)$ be an algebraic quantum hypergroup.  If furthermore $\Delta$ is an algebra homomorphism, then $(A,\Delta)$ is an algebraic quantum group.
\snl\bf Proof: \rm 
So, as before, let $(A,\Delta)$ be an algebraic quantum hypergroup with counit $\varepsilon$, a left integral $\varphi$ and an antipode $S$, relative to $\varphi$. Moreover, assume that $\Delta$ is an algebra homomorphism. We claim that  
$$\align m((S \otimes \iota)(\Delta(x) (1 \otimes y))) &= \varepsilon (x) y \\
m((\iota \otimes S)((x \otimes 1) \Delta(y))) &= \varepsilon(y) x
\endalign$$
for all $x, y \in A$ where $m : A \otimes A \rightarrow A$ is the multiplication map. We will only prove the first formula because the proof of the second one is completely similar.  
\snl
As before, define $\psi = \varphi \circ S$, so that $\psi$ is a faithful right integral on $A$. Take $a,b\in A$ and put $x = (\psi \otimes \iota)((b \otimes 1)\Delta(a))$ for $a, b \in A$.  Then we have for all $y \in A$ that
$$\align
\Delta(x)(1 \otimes y) &= (\psi \otimes \iota \otimes \iota)((\iota \otimes \Delta)((b \otimes 1)\Delta(a)) (1 \otimes 1 \otimes y))\\
&= (\psi \otimes \iota \otimes \iota)((b \otimes 1 \otimes 1) (\Delta \otimes \iota) (\Delta(a)(1 \otimes y))).
\endalign$$
Apply $S \otimes \iota$ on both sides of this equation.  Using Proposition 2.2, we obtain
$$(S \otimes \iota)(\Delta(x)(1 \otimes y)) = (\psi \otimes \iota \otimes \iota)(\Delta_{12}(b) \Delta_{13}(a) (1 \otimes 1 \otimes y)).$$
Therefore we have,
$$\align 
m(S \otimes \iota) (\Delta(x) (1 \otimes y)) &= (\psi \otimes \iota)(\Delta(b) \Delta(a) (1 \otimes y))\\
&= (\psi \otimes \iota)(\Delta(ba) (1 \otimes y)) = \psi(ba)y = \varepsilon(x) y.
\endalign$$
Now $(A,\Delta)$ is an algebra with an ordinary regular comultiplication in the sense of [VD2]. By the use of [VD2, Proposition 2.9], we obtain that $(A,\Delta)$ is a regular multiplier Hopf algebra (with integrals). \hfill$\blacksquare$
\einspr

We have seen already in Section 1 (Proposition 1.14) that any algebraic quantum group is also an algebraic quantum hypergroup. If we combine this result with the one in the previous proposition, we see that {\it the algebraic quantum groups are precisely those algebraic quantum hypergroups with a comultiplication that is an algebra homomorphism}. A similar statement is true for $^*$-algebraic quantum hypergroups.
\nl
Next, we prove the following uniqueness result for algebraic quantum hypergroups. Again, $(A,\Delta)$ is an algebraic quantum hypergroup with counit $\varepsilon$, left integral $\varphi$ and antipode $S$, relative to $\varphi$. 

\inspr{2.4} Proposition \rm 
If $\varphi'$ is another left invariant functional on $(A, \Delta)$, then $\varphi' = \lambda \varphi$ for some scalar $\lambda \in \Bbb C$.
\snl\bf Proof: \rm 
Take $a, b \in A$ and apply $\varphi'$ to both expressions in the equation
$$S(\iota \otimes \varphi) (\Delta (a) (1 \otimes b)) = (\iota \otimes \varphi) ((1 \otimes a) \Delta(b))$$
for all $a, b \in A$. Because $\varphi' \circ S$ is right invariant, the left hand side will give $\varphi'(S(a)) \varphi(b)$.  For the right hand side we get $\varphi (a \delta_b)$ where $\delta_b$ is defined in $M(A)$ by the formula $\delta_b = (\varphi' \otimes \iota) \Delta(b)$. Because $\varphi'(S(a)) \varphi(b) = \varphi(a\delta_b)$ for all $a \in A$ and because $\varphi$ is faithful, we must have a multiplier $\delta \in M(A)$ such that $\delta_b = \varphi(b) \delta$ for all $b\in A$.  If we apply $\varepsilon$ we get $\varphi'(b) = \varphi(b) \varepsilon(\delta)$ for all $b \in A$ and with $\lambda = \varepsilon (\delta)$, we find the desired result.\hfill$\blacksquare$
\einspr

Proposition 2.4 not only proves the uniqueness of the left integral and the right integral (by composing with $S$), it also proves the uniqueness of the antipode in the following sense.  If $\varphi$ and $\varphi'$ are faithful left integrals and $S$ and $S'$ antipodes relative to $\varphi$ and $\varphi'$ respectively, then the above result gives that $\varphi'$ is a scalar multiple of $\varphi$ and also that $S'$ must be the same as $S$, because the antipode is uniquely determined by the faithful left integral.
\snl
So, given the algebraic quantum hypergroup $(A,\Delta)$, we have a uniquely defined counit $\varepsilon$ and antipode $S$, as well as (up to a scalar) a unique left integral $\varphi$ and a unique right integral $\psi$. This justifies the way we defined algebraic quantum hypergroups in Definition 1.10. 
\snl
While proving the uniqueness in Proposition 2.4, we also have shown already part of the following result.

\inspr{2.5} Proposition \rm 
There is a unique invertible element $\delta \in M(A)$ such that 
\snl
(1) $(\varphi \otimes \iota) \Delta(a) = \varphi(a)\delta$ and $(\iota \otimes \psi) \Delta(a) = \psi(a) \delta^{-1}$ \newline 
(2) $\varphi(S(a)) = \varphi(a\delta)$ 
\snl
for all $a\in A$. We also have $\varepsilon(\delta) = 1$ and $S(\delta) = \delta^{-1}$.
\snl\bf Proof: \rm
If in the proof of Proposition 2.4, we take $\varphi' = \varphi$ we get a multiplier $\delta \in M(A)$ such that $(\varphi \otimes \iota) \Delta(a) = \varphi(a)\delta$ and $\varphi(S(a)) = \varphi(a\delta)$ for all $a \in A$.  This gives the first part of (1) and (2).  If we apply $\varepsilon$ on the first equation, we find $\varepsilon(\delta) = 1$. Because $S$ flips the coproduct and if we let $\psi = \varphi \circ S$, we get $(\iota \otimes \psi)\Delta(a) = \psi (a) \delta'$ where $\delta' = S^{-1}(\delta)$.  So, it remains to show that $\delta' = \delta^{-1}$.
\snl
If we apply $\varphi$ to the formula in Proposition 2.2, we get $\varphi(S(a)) \psi(b\delta') = \varphi(b) \psi(a)$ for all $a, b \in A$. In particular, we have $\varphi(b) = \psi(b\delta')$ for all $b \in A$.
\snl
Therefore, we have $\varphi(b) = \varphi(S(b\delta')) = \varphi(b\delta'\delta)$ and so $\delta'\delta = 1$.  On the other hand, we also have $\psi(b) = \varphi (S(b)) = \varphi(b\delta) = \psi(b\delta\delta')$ and so we have $\delta\delta' = 1$.  Therefore $\delta$ is invertible and $\delta^{-1} = \delta' = S^{-1}(\delta)$, or equivalently $S(\delta) = \delta^{-1}$.\hfill $\blacksquare$
\einspr

In the $^*$-algebra case, we have seen that we can always assume that $\varphi$ is self-adjoint. This will imply that $\delta^*=\delta$. 
\snl
The multiplier $\delta$ is called the {\it modular element}. The terminology comes from the theory of locally compact groups where the function relating the left with the right Haar measure is called the modular function.
\nl
Because $S$ is an anti-isomorphism that flips the coproduct, the square $S^2$ of the antipode will be an isomorphism that leaves the coproduct invariant. It follows that the composition $\varphi\circ S^2$ of the left integral $\varphi$ with $S^2$ will again be a left integral. By the uniquess of left integrals, we must have a complex number $\tau$ satisfying $\varphi(S^2(a))=\tau \varphi(a)$ for all $a\in A$. This number is called the {\it scaling constant}. In the $^*$-algebra case, one can show that $|\tau|=1$.
\nl
Finally, just as in the case of algebraic quantum groups, also here we have the existence of the modular automorphisms, as in the following proposition.

\inspr{2.6} Proposition \rm 
There is a unique automorphism $\sigma$ of $A$ such that $\varphi(ab) = \varphi(b\sigma(a))$ for all $a, b \in A$. We also have $\varphi(\sigma(a))=\varphi(a)$ for all $a$ in $A$. Similarly, there is a unique automorphism $\sigma'$ of $A$ satisfying $\psi(ab)=\psi(b\sigma'(a))$ for all $a,b\in A$. Also here $\psi(\sigma'(a))=\psi(a)$ for all $a$. 
\snl \bf Proof: \rm  
For all $p$, $q$, $x$ in $A$ we have
$$\align
(\psi \otimes \varphi) ((1 \otimes p) (\iota \otimes S) ((x \otimes 1)\Delta(q))) 
   &= \varphi(p(\psi \otimes \iota) (\iota \otimes S) ((x \otimes 1)\Delta (q)))\\
   &= \varphi(p(\psi \otimes \iota) (\Delta(x) (q \otimes 1)))\\
   &= \psi ((\iota \otimes \varphi) ((1 \otimes p) \Delta(x) (q \otimes 1)))\\
   &= \psi (S((\iota \otimes \varphi)(\Delta(p) (1 \otimes x))) q) \\
   &= \varphi((((\psi \circ S) \otimes \iota) ((S^{-1} (q) \otimes 1)\Delta(p)))x).
\endalign$$
On the other hand, we also have
$$(\psi \otimes \varphi) ((1 \otimes p) (\iota \otimes S)((x \otimes 1) \Delta(q))) = \psi(x(\iota \otimes (\varphi \circ S)) (\Delta(q) (1 \otimes S^{-1} (p)))).$$
Now assume that $\psi = \varphi \circ S$. Then we get $\psi \circ S = \tau \varphi$ as well as $\psi(y)=\varphi(y\delta)$.
Then the above calculation will give us the formula $\varphi(ax) = \varphi(xb)$ for all $x\in A$ where $a = (\varphi \otimes \iota)((S^{-1} (q) \otimes 1) \Delta(p))$ and $b = \frac{1}{\tau} (\iota \otimes \psi)(\Delta(q) (1 \otimes S^{-1} (p)))\delta$.
\snl
Because $\varphi$ is supposed to be faithful, the element $b$ is uniquely determined by the element $a$ and therefore, we can define $\sigma(a)=b$. Moreover, because all elements in $A$ are of the form $a$ above, the map $\sigma$ is defined on all of $A$. The map $\sigma$ is injective, again by the faithfulness of $\varphi$. It is also surjective because all elements in $A$ are also of the form $b$ above.
\snl
To show that $\sigma$ is a homomorphism, take elements $a$, $b$ and $c$ in $A$. We have $\varphi(abc) = \varphi(a(bc)) = \varphi((bc) \sigma(a)) = \varphi (b(c\sigma(a))) = \varphi ((c \sigma(a)) \sigma(b)) = \varphi (c(\sigma(a) \sigma(b)))$.
Because also $\varphi(abc) = \varphi((ab)c) = \varphi(c\sigma(ab))$, it follows from the faithfulness of $\varphi$, that $\sigma(ab) = \sigma(a) \sigma(b)$ and so $\sigma$ is an algebra homomorphism. 
\snl
When we apply the result two times, we get
$$\varphi(ab)=\varphi(b\sigma(a))=\varphi(\sigma(a)\sigma(b))=\varphi(\sigma(ab))$$
for all $a,b$ in $A$. Because $A^2=A$, as a consequence of Proposition 1.6, we get that $\varphi$ is $\sigma$-invariant. 
\snl
This proves the statement about $\varphi$ and $\sigma$. Using that $\psi=\varphi\circ S$ one gets easily the statement for $\psi$, taking $\sigma'=S^{-1}\circ\sigma^{-1}\circ S$.
\hfill $\blacksquare$
\einspr

The automorphisms $\sigma$ and $\sigma'$ are called the {\it modular automorphisms} of $A$ associated with $\varphi$ and $\psi$ respectively. The terminology comes from the theory of operator algebras. In general, the algebra $A$ is not abelian and the integrals are not traces. The modular automorphisms take care of the possible problems that arise from these facts. In particular, the results of the above proposition will be necessary for proving elementary properties of the dual algebraic quantum hypergroup as we will introduce in the next section.
\snl
There are various extra properties one can easily deduce from the above propositions. We collect all of them in the next proposition. We also prove several relations with the other data associated with an algebraic quantum hypergroup.

\inspr{2.7} Proposition \rm With the notations of before, we get:
\snl
(1) The modular automorphisms $\sigma$ and $\sigma'$ are related with each other in two ways. \newline
$\text{ }$\hskip 0.5cm We have $\sigma\circ S \circ\sigma' = S$ but also $\sigma'(a)=\delta \sigma(a)\delta^{-1}$. \newline
(2) We have  $\sigma(\delta)=\frac{1}{\tau}\delta$ as well as $\sigma'(\delta)=\frac{1}{\tau}\delta$. \newline
(3) The modular automorphisms $\sigma$ and $\sigma'$ commute with each other. \newline
(4) The modular automorphisms $\sigma$ and $\sigma'$ commute with $S^2$. \newline 
(5) For all $a$ we have $\Delta(\sigma(a)) = (S^2 \otimes \sigma) \Delta(a)$ and $\Delta(\sigma'(a)) = (\sigma' \otimes S^{-2}) \Delta (a)$. \newline
(6) For all $a$ we have also $\Delta(S^2(a)) = (\sigma \otimes {\sigma'}^{-1}) \Delta(a)$.

\snl \bf Proof: \rm
The first statement in (1) is already found in the proof of the previous proposition. We used this formula to define $\sigma'$. The second formula in (1) is obtained from the fact that $\psi(a)=\varphi(a\delta)$.
\snl
To prove (2) we use that
$$\varphi(S^2(a)) = \varphi(S(a)\delta) = \varphi(S(\delta^{-1}a)) = \varphi(\delta^{-1} a \delta)$$
for all $a$ in $A$ and we find $\tau\varphi(a)= \varphi(\delta^{-1} a \delta)$. This implies that $\sigma(\delta)=\frac{1}{\tau}\delta$. From the relations between $\sigma$ and $\sigma'$ from (1), and using that $S(\delta)=\delta^{-1}$ we also get $\sigma'(\delta)=\frac{1}{\tau}\delta$.
\snl
Using some of the results above we find
$$\sigma(\sigma'(a)) = \sigma(\delta \sigma(a)\delta^{-1}) = \sigma(\delta) \sigma^2(a) \sigma(\delta^{-1}) = \delta \sigma^2 (a) \delta^{-1} = \sigma' (\sigma(a))$$
and so the automorphisms $\sigma$ and $\sigma'$ commute with each other.
\snl
Statement (4) can be shown in different ways. If we combine e.g. the results in (1) and (2) with each other, we find that not only $\sigma\circ S= S\circ {\sigma'}^{-1}$ but also that  ${\sigma'}^{-1}\circ S=S\circ \sigma$. The two together give us that $\sigma$ and $S^2$ commute. 
\snl
We now prove the formula in (5) for $\sigma$.  The other one follows e.g.\ from the first relation between $\sigma$ and $\sigma'$ in (1) and the fact that $S$ flips the antipode.  For all $a$, $b$ in $A$, we have
$$\align
(\iota \otimes \varphi) ((1 \otimes b) (S^2 \otimes \sigma) \Delta(a)) 
  &= S^2 (\iota \otimes \varphi)((1 \otimes b) (\iota \otimes \sigma) \Delta(a))\\
  &= S^2 (\iota \otimes \varphi) (\Delta(a)(1 \otimes b)) \\
  &= S (\iota \otimes \varphi)((1 \otimes a) \Delta(b))\\
  &= S (\iota \otimes \varphi)(\Delta(b) (1 \otimes \sigma(a)))\\ 
  &= (\iota \otimes \varphi) ((1 \otimes b) \Delta(\sigma(a))).
\endalign$$
Now the result follows from the faithfulness of $\varphi$ on $A$.
\snl
(6)  If we aply $\varepsilon$ to the second leg in the formula for $\Delta\circ\sigma$ and to the first leg in the formula for $\Delta\circ\sigma'$ in (5) we get 
$$\align S^{-2}\sigma(a)&=(\iota\ot (\varepsilon\circ\sigma))\Delta(a)\\
  S^2\sigma'(a)&=((\varepsilon\circ\sigma') \ot \iota)\Delta(a)
\endalign$$
for all $a$. From $\varepsilon(\delta)=1$ and  $\sigma'(a)=\delta\sigma(a)\delta^{-1}$ it follows that $\varepsilon\circ\sigma=\varepsilon\circ\sigma'$. Therefore 
$$(\iota\ot(\varepsilon\circ\sigma) \ot \iota)\Delta^{(2)}(a)=(\iota\ot(\varepsilon\circ\sigma') \ot \iota)\Delta^{(2)}(a)$$
and if we use the two previous formulas, we get from this that 
$$(S^{-2}\sigma \ot\iota)\Delta(a)=(\iota\ot S^2\sigma')\Delta(a).$$
Because $\Delta(S^2(a))=(S^2 \ot S^2)\Delta(a)$, we get the desired formula.
\hfill$\blacksquare$
\einspr

There is also the formula $\Delta(\delta)=\delta\otimes\delta$ but this is more subtle. Because $\Delta$ is no longer assumed to be an algebra homomorphism, it is not at all obvious how to possibly extend it to the multiplier algebra $M(A)$. We will come back to this problem after we have obtained duality. Then also the formula for $\Delta(\delta)$ will be considered (see a remark following Proposition 4.1).
\snl
Remark that the last property in the above proposition, in the case of ordinary algebraic quantum groups, was first proven in [K-VD], using among other things that $\Delta(\delta)=\delta\otimes\delta$. The proof there is also more complicated. The proof we give here is the same as the one for algebraic quantum groups given in the appendix of [L-VD1].
\nl
It is not useful to look at the motivating example to illustrate these results. Indeed, in Example 1.11, the algebra is abelian so that the modular automorphisms $\sigma$ and $\sigma'$ are trivial. Also $\varphi=\psi$ so that $\delta=1$. Finally $\tau=1$ because $S^2=\iota$. We need more complicated examples to illustrate these results (see [D-VD2]).
\snl
If we have a $^*$-algebraic quantum group of compact type with a positive left integral $\varphi$, then $\varphi(1)>0$ and a standard argument gives that $\psi=\varphi$. In particular $\sigma=\sigma'$ and $\delta=1$. Further, the objects may be non-trivial  as this is already the case for $^*$-algebraic quantum groups of compact type with positive integrals (i.e.\ compact quantum groups in the sense of [W1] and [W2]).  

\nl\nl

\bf 3. Duality and biduality \rm
\nl
Let $(A, \Delta)$ be an algebraic quantum hypergroup (in the sense of Definition 1.10). In this section, we will construct the dual $(\widehat A, \widehat\Delta)$ and we will show that it is again an algebraic quantum hypergroup. The construction goes very much  as in the case of ordinary algebraic quantum groups (cf.\ [VD2]). So, also here we start with defining the following subspace of the dual space $A'$.

\inspr{3.1} Definition \rm
Let $\varphi$ be a left integral on $(A,\Delta)$. We define $\widehat A$ as the space of linear functionals on $A$ of the form $\varphi(\,\cdot\,a)$ where $a\in A$.
\einspr

Because of the Propositions 2.5 and 2.6, we get
 $$\{\varphi(a\,\cdot\,)\mid a \in A\} = \{\varphi(\,\cdot\, a)\mid a \in A\} = \{\psi (\,\cdot\, a) \mid a \in A\} = \{\psi(a \, \cdot\,)\mid a \in A\}$$
where $\varphi$ is a left integral and $\psi$ is a right integral on $A$. Therefore, any element in $\widehat A$ can be written in any of the four different ways above. We will use freely any of these expressions as, depending on the case, each of them is useful. Occasionally, we will use $\langle \omega, x\rangle$ to denote the value $\omega(x)$ of an element $\omega\in \widehat A$ (or even $\omega\in A'$) in an element $x\in A$. By the faithfulness of the integrals, the space $\widehat A$ is separating which implies that the pairing between $A$ and $\widehat A$ is non-degenerate.
\snl
We will prove that $\widehat{A}$ can again be made into an algebraic quantum hypergroup. Furthermore, by considering the dual of $\widehat{A}$, i.e.\ the bidual of $A$, we will recover the original algebraic quantum hypergroup $A$. 
\snl
We start by making $\widehat A$ into an algebra by dualizing the coproduct. This is done in the following proposition.

\inspr{3.2} Proposition \rm 
For $\omega, \omega' \in \widehat{A}$, we can define a linear functional $\omega\omega'$ on $A$ by the formula $(\omega \omega')(x) = (\omega \otimes \omega') \Delta(x)$ for all $x \in A$.  We get $\omega\omega' \in \widehat{A}$. This product on $\widehat A$ is associative and non-degenerate.

\snl\bf Proof: \rm Let $\omega,\omega'\in \widehat A$ and assume that  $\omega' = \varphi(\,\cdot\, a)$ with $a\in A$. We have
$$\align
(\omega \omega')(x) &= (\omega \otimes \varphi(\,\cdot\, a)) (\Delta(x)) = (\omega \otimes \varphi)(\Delta(x) (1 \otimes a))\\
&= (\omega \circ S^{-1}) ((\iota \otimes \varphi)((1 \otimes x) \Delta(a))) = \varphi(x((\omega \circ S^{-1}) \otimes \iota) \Delta(a)).
\endalign$$
We see that not only the product $\omega\omega'$ is well-defined as a linear functional on $A$, but also that it has the form $\varphi(\,\cdot\,b)$ with $b=((\omega \circ S^{-1}) \otimes \iota) \Delta(a)$. Because also $\omega \in \widehat{A}$, we have $b\in A$ (and not just in $M(A)$).  So $\omega \omega'\in \widehat{A}$. Therefore, we have defined a product in $\widehat A$.
\snl
The associativity of this product in $\widehat{A}$ is an easy consequence of the coassociativity of $\Delta$ on $A$.  To prove that the product is non-degenerate, assume that $a\in A$ and that $\omega \varphi(\,\cdot\, a) = 0$ for all $\omega \in \widehat{A}$.  This implies $\omega(S^{-1} (a)) = 0$ for all $\omega \in \widehat{A}$.  As $\widehat A$ is separating points of $A$, we conclude that $a = 0$.  Similarly, $\omega\varphi(\,\cdot\, a) = 0$ for all $a \in A$ implies that $\omega = 0$. \hfill $\blacksquare$
\einspr

If $A$ is a $^\ast$-algebra, we assume that $\Delta(a^\ast) = \Delta(a)^\ast$ for all $a \in A$. In this case, we can also define an involution on $\widehat A$. We let $\omega^*(x)=\omega(S(x)^*)^-$ when $x\in A$ and $\omega\in \widehat A$ (where $\lambda^-$ is the complex conjugate of the complex number $\lambda$). It is not hard to see that $\omega^*$ is again in $\widehat A$ and that we have made $\widehat A$ into an involutive algebra. Among other properties, one has to use that $S(S(x)^\ast)^\ast = x$ for all $x \in A$.
\snl
We have seen that the elements of $\widehat{A}$ can be expressed in four different forms.  When we use these different forms in the definition of the product in $\widehat{A}$, we get the following useful expressions. As before, $\varphi$ is a left integral and $\psi$ is a right integral.

\inspr{3.3} Proposition \rm
Whenever we have $a \in A$ and $\omega \in \widehat{A}$, we get
$$\alignat{2}
(1)\qquad &\omega \varphi(\,\cdot\, a) = \varphi(\,\cdot\, b)& \qquad \text{with} \qquad & b = ((\omega \circ S^{-1}) \otimes \iota) \Delta(a)\\
(2)\qquad  &\omega \varphi(a \,\cdot\,) = \varphi(c\,\cdot\,) & \qquad \text{with} \qquad & c= ((\omega \circ S) \otimes \iota) \Delta(a)\\
(3)\qquad  &\psi(\,\cdot\, a) \omega = \psi(\,\cdot\, d) & \qquad \text{with} \qquad & d = (\iota \otimes (\omega \circ S)) \Delta(a)\\
(4)\qquad  &\psi (a \,\cdot\,) \omega = \psi (e\,\cdot\,) & \qquad \text {with} \qquad & e = (\iota \otimes (\omega \circ S^{-1})) \Delta(a)
\endalignat$$

\snl\bf Proof: \rm
The first formula was already obtained in the proof of the previous proposition. The 3 other formulas are obtained in a completely similar fashion, now using not only the definition of $S$, involving $\varphi$ (cf. Definition 1.9), but also the other relation, involving $\psi$ (see Proposition 2.2). \hfill $\blacksquare$
\einspr

The above formulas can also be used to get products $f\omega$ and $\omega f$, for all $f\in A'$ and $\omega \in \widehat{A}$. We can e.g.\ define $f\omega$ when $\omega=\varphi(\,\cdot\, a)$ with $a\in A$ (using the alternative notation) by
$$
 \langle f \omega, x \rangle = \langle f \otimes \varphi(\,\cdot\, a), \Delta (x)\rangle = \langle f \otimes \varphi, \Delta (x) (1 \otimes a)\rangle$$
and we see (with the argument as in the proof of Proposition 3.2) that indeed, this is $\varphi(xb)$ with $b=((f \circ S^{-1}) \otimes \iota) \Delta(a)$.  In general, these products will not belong to $\widehat A$ anymore. However, we clearly have that $A'$ is a $\widehat A$-bimodule.
\snl
On the other hand, we have the following result.

\inspr{3.4} Proposition \rm 
Let $f$ in $A'$ such that $(f\ot\iota)\Delta(a)$ and $(\iota\ot f)\Delta(a)$ belong to $A$ (and not only to $M(A)$) for all $a$ in $A$. Then $f\omega$ and $\omega f$ belong to $\widehat A$ for all $\omega\in \widehat A$. This defines an element in $M(\widehat A)$. All of $M(\widehat A)$ can be realized in this way. 
\snl \bf Proof: \rm
From the argument above, we see that $f\omega=\varphi(\,\cdot\,b)$ when $\omega=\varphi(\,\cdot\,a)$ with $a\in A$ and $b=((f \circ S^{-1}) \otimes \iota) \Delta(a)$. Using Proposition 2.1 it follows that $b\in A$ when $f$ satisfies the condition in the formulation of the proposition. Therefore, $f\omega\in\widehat A$ when $\omega\in\widehat A$. Similarly for $\omega f$. So, we see that functionals like $f$ give rise to multipliers in $M(\widehat A)$.
\snl
Conversely, suppose that $m$ is a multiplier in $M(\widehat A)$ and assume that $m\omega=\varphi(\,\cdot\,b)$ for $\omega=\varphi(\,\cdot\, a)$ with $a,b\in A$. Define a linear functional $f$ on $A$ by $f(S^{-1}(a))=\varepsilon(b)$. From the fact that $m(\omega\omega')=(m\omega)\omega'$ for all $\omega,\omega'\in \widehat A$, it can be shown that indeed $m\omega=f\omega$ for all $\omega\in\widehat A$. It also follows easily that $f$ satisfies the conditions of the proposition.\hfill $\blacksquare$
\einspr

As before, we see that $f\omega=0$ for all $\omega\in \widehat A$ will imply $f=0$. So, the above map can be used to identify $M(\widehat A)$ with the subspace of elements $f$ in $A'$ with the given property. As a consequence, the pairing between $A$ and $\widehat A$ can be extended to a pairing between $A$ and $M(\widehat A)$ (in the sense of bilinear maps). Essentially by definition, we get
$$\align \langle f\omega, x\rangle &= \langle f, (\iota\otimes\omega)\Delta(x)\rangle \\
          \langle \omega f, x\rangle &= \langle f, (\omega\otimes\iota)\Delta(x)\rangle
\endalign$$
for all such functionals $f$ and all $x\in A$. If we consider the extended pairing, we can consider $f$ as an element of $M(\widehat A)$ in these formulas. We will come back to these formulas when we consider the module structures in the next section. 
\snl
Observe that the counit $\varepsilon$, as a linear functional on $A$, is in fact the unit in the multiplier algebra $M(\widehat A)$ of $\widehat A$. This follows from the formulas $(\varepsilon\otimes\iota)\Delta(a)=a$ and $(\iota\otimes\varepsilon)\Delta(a)=a$ for all $a\in A$.
\snl 
In a similar way, we can consider elements in $M(\widehat A\ot\widehat A)$ as linear functionals on $A\ot A$. This will be helpful to understand the coproduct on $\widehat A$.
\nl
Let us now define this {\it comultiplication} $\widehat{\Delta}$ on $\widehat{A}$. Roughly speaking, and when considering elements in $M(\widehat A\otimes\widehat A)$ as linear functionals on $A\ot A$, the coproduct is dual to the multiplication in $A$ in the sense that $\langle\widehat\Delta(\omega),x\ot y\rangle=\langle \omega, xy \rangle$ when $x,y\in A$. However, because we have that $\widehat\Delta$ does not map into $\widehat A\ot \widehat A$ but rather into the multiplier algebra of this tensor product, we have to be more careful. We will define the coproduct by giving the expressions for $(\omega_1 \otimes 1)\widehat{\Delta} (\omega_2)$ and $\widehat{\Delta}(\omega_1)(1 \otimes \omega_2)$ for all $\omega_1$, $\omega_2$ in $\widehat{A}$. These objects will be in $\widehat{A}\otimes \widehat{A}$. 

\inspr{3.5} Definition \rm 
Let $\omega_1,\omega_2 \in \widehat{A}$.  Then we put
$$\align
\langle(\omega_1 \otimes 1)\widehat{\Delta}(\omega_2), x \otimes y\rangle &= \langle\omega_1 \otimes \omega_2, \Delta(x) (1 \otimes y)\rangle\\
\langle\widehat{\Delta}(\omega_1) (1 \otimes \omega_2), x \otimes y\rangle &= \langle\omega_1 \otimes \omega_2,(x \otimes 1) \Delta(y)\rangle
\endalign$$
for all $x, y \in A$.
\einspr

We will first show  that the functionals in Definition 3.5 are well-defined and again in $\widehat{A} \otimes \widehat{A}$. Then, it will be possible to define $\widehat \Delta(\omega)$ in $M(\widehat A\ot \widehat A)$ and we will get the expected formula. 

\inspr{3.6} Lemma \rm 
If $\omega_1, \omega_2 \in \widehat{A}$, then $(\omega_1 \otimes 1) \widehat{\Delta}(\omega_2)$ and $\widehat{\Delta}(\omega_1)(1 \otimes \omega_2)$ (as defined in 3.5 above) are in $\widehat{A} \otimes \widehat{A}$.  These two formulas define $\widehat{\Delta}(\omega)$ as a multiplier in $M(\widehat{A}\otimes \widehat{A})$ for all $\omega \in \widehat{A}$. 

\snl\bf Proof: \rm 
Let $\omega_1 = \psi (a \,\cdot\,)$ and $\omega_2 = \psi(b\,\cdot\,)$ where $a,b \in A$.  For all $x, y \in A$, we have
$$\align
\langle(\omega_1 \otimes 1)\widehat{\Delta}(\omega_2), x \otimes y \rangle &= (\omega_1 \otimes \omega_2) (\Delta(x) (1 \otimes y))\\
 &= \omega_2 (\psi \otimes \iota) ((a \otimes 1) \Delta(x) (1 \otimes y)) \\
 &= \omega_2 (S^{-1} ((\psi \otimes \iota)(\Delta(a) (x \otimes 1))) y)\\
&= \psi ((\psi \otimes \iota) ((1 \otimes b)(\iota \otimes S^{-1}) (\Delta(a) (x \otimes 1))) y)\\
&= \psi ((\psi \otimes \iota) ((\iota \otimes S^{-1}) (\Delta(a) (x \otimes S(b)))) y)\\
&= \psi ((\psi \otimes \iota) ((\iota \otimes S^{-1}) (\Delta(a) (x \otimes S(b))) (1 \otimes y)))\\
&= (\psi \otimes \psi)((\iota \otimes S^{-1})(\Delta(a) (1 \otimes S(b))) (x \otimes y))\\
\endalign$$
We obtain that $(\omega_1 \otimes 1) \widehat{\Delta} (\omega_2)$ is a well-defined element in $\widehat{A} \otimes \widehat{A}$.
\snl
To prove that also $\widehat{\Delta} (\omega_1) (1 \otimes \omega_2)$ is well defined in $\widehat{A} \otimes \widehat{A}$, we use the expressions $\omega_1 = \varphi(\,\cdot\, a)$ and $\omega_2 = \varphi (\,\cdot\, b)$ where $\varphi$ is a left integral on $A$.  Then we find that 
$$\langle\widehat{\Delta} (\omega_1) (1 \otimes \omega_2),x\ot y \rangle = (\varphi \otimes \varphi)((x\ot y)(S^{-1} \otimes \iota)((S(a) \otimes 1) \Delta(b)))$$
for all $x,y \in A$.
\snl
Using that the product in $\widehat A$ is dual to the coproduct on $A$, it easily follows from the definitions above that   $((\omega_1 \otimes 1)\widehat{\Delta}(\omega_2))(1 \otimes \omega_3) = (\omega_1 \otimes 1)(\widehat{\Delta}(\omega_2)(1 \otimes \omega_3))$  for all $\omega_1$, $\omega_2$ and $\omega_3$ in $\widehat{A}$.  Therefore, $\widehat{\Delta}(\omega)$ is defined as a two-sided multiplier in $M(\widehat{A} \otimes \widehat{A})$ for all $\omega \in \widehat{A}$. \hfill $\blacksquare$
\einspr

Let us now argue that $\langle\widehat{\Delta}(\omega),x \otimes y\rangle = \langle\omega,xy \rangle$ for all $x, y \in A$. Consider the first formula in the definition above. Because of one of the remarks, following Proposition 3.4, we have
$$\langle(\omega_1\otimes 1)\widehat\Delta(\omega_2),x'\otimes y \rangle=
                    \langle\widehat \Delta(\omega_2),(\omega_1\otimes 1)\Delta(x') \otimes y \rangle$$
and so
$$\align \langle\widehat\Delta(\omega_2),(\omega_1\otimes\iota)\Delta(x') \otimes y\rangle 
                  &=\langle \omega_1 \otimes \omega_2, \Delta(x')(1\otimes y)\rangle \\
                  &=\omega_2(((\omega_1\otimes\iota)\Delta(x'))y)
\endalign$$
whenever $\omega_1,\omega_2\in \widehat A$ and $x',y\in A$. Therefore, we get the desired formula for elements $x,y\in A$ with $x$ of the form $(\omega_1\otimes \iota)\Delta(x')$. However, all elements in $A$ are of this form and so we have shown the formula for all pairs. 

\inspr{3.7} Proposition \rm 
The map $\widehat{\Delta}:\widehat{A} \rightarrow M(\widehat{A} \otimes \widehat{A})$ is a regular comultiplication on $\widehat{A}$. 

\snl\bf Proof: \rm 
For all $\omega_1$, $\omega_2$ and $\omega_3$ in $\widehat{A}$ and $x$, $y$, $z$ in $A$, we have
$$\align
\langle(\omega_1 \otimes 1 \otimes 1)(\widehat{\Delta} \otimes \iota)(\widehat{\Delta} (\omega_2)
&(1 \otimes \omega_3)), x \otimes y \otimes z\rangle \\
&= \langle \omega_1 \otimes (\widehat{\Delta} (\omega_2) (1 \otimes \omega_3)), (\Delta(x) (1 \otimes y)) \otimes z\rangle\\
&= \langle \omega_1 \otimes \omega_2 \otimes \omega_3, \Delta_{12} (x) (1 \otimes y \otimes 1) \Delta_{23} (z)\rangle\\
&= \langle ((\omega_1 \otimes 1)\widehat{\Delta}(\omega_2)) \otimes \omega_3, x \otimes ((y \otimes 1) \Delta(z))\rangle\\
&= \langle (\iota \otimes \widehat{\Delta})((\omega_1 \otimes 1)\widehat{\Delta}(\omega_2)) (1 \otimes 1 \otimes \omega_3), x \otimes y \otimes z\rangle.\endalign$$
This shows that $\widehat{\Delta}$ is coassociative in the sense of Definition 1.1. Remark that we have used the 'leg-numbering' notation as explained earlier (in the proof of Proposition 2.1). 
\snl
To show that this coproduct is regular, we prove that also the elements $\widehat{\Delta} (\omega_1) (\omega_2 \otimes 1)$ and $(1 \otimes \omega_1) \widehat{\Delta}(\omega_2)$ are in $\widehat{A} \otimes \widehat{A}$ for all $\omega_1$ and $\omega_2$ in $\widehat{A}$. With $x,y\in A$, we have
$$\align
\langle \widehat{\Delta}(\omega_1) (\omega_2 \otimes 1), x \otimes y \rangle &= \langle \widehat{\Delta}(\omega_1), (\iota \otimes \omega_2) \Delta(x) \otimes y\rangle\\
&= \langle \omega_1, (\iota \otimes \omega_2)(\Delta(x) (y \otimes 1))\rangle \\
&= \langle \omega_1 \otimes \omega_2, \Delta(x) (y \otimes 1)\rangle.
\endalign$$
Observe that in the argument above, we have considered $\widehat\Delta(\omega_1)$ as a linear functional on $A\ot A$ as above. If now $\omega_1 = \varphi(a \,\cdot\,)$ and $\omega_2 = \varphi(\,\cdot\, b)$ with $a,b \in A$, we have
$$\align
\langle \widehat{\Delta} (\omega_1)(\omega_2 \otimes 1), x \otimes y\rangle
&= \omega_1 (((\iota \otimes \varphi) (\Delta(x) (1 \otimes b))) y)\\
&= \omega_1 ((S^{-1} (\iota \otimes \varphi)((1 \otimes x) \Delta(b)))y)\\
&= \varphi(((\iota \otimes \varphi)(1 \otimes x) (S^{-1} \otimes \iota)(\Delta(b) (S(a)\otimes 1)))y)\\
&= \langle \varphi \otimes \varphi, (1 \otimes x) (S^{-1} \otimes \iota) (\Delta(b) (S(a) \otimes 1)) (y \otimes 1)\rangle.
\endalign$$
Write $(S^{-1} \otimes \iota) (\Delta(b) (S(a) \otimes 1)) = \sum p_i \otimes q_i$ in $A \otimes A$. Then we have shown that $\widehat{\Delta}(\omega_1)(\omega_2 \otimes 1) = \sum_i \varphi(\,\cdot\, q_i) \otimes \varphi(p_i \,\cdot\,)$. This proves the first claim. To prove the second statement, we use 
$$\langle (1 \otimes \omega_1)\widehat{\Delta}(\omega_2), x \otimes y\rangle = 
\langle \omega_1 \otimes \omega_2, (1 \otimes x) \Delta(y)\rangle.$$
Now we take $\omega_1 = \psi(a \,\cdot\,)$ and $\omega_2 = \psi (\,\cdot\, b)$ with $a,b\in A$.  Then we  obtain $(1 \otimes \omega_1)\widehat{\Delta} (\omega_2) = \sum_i \psi (\,\cdot\, q_i) \otimes \psi (p_i \,\cdot\,)$ where $\sum_i p_i \otimes q_i = (\iota \otimes S^{-1}) ((1 \otimes S(b)) \Delta(a))$. This completes the proof. \hfill$\blacksquare$
\einspr

If $A$ is a $^\ast$-algebraic quantum hypergroup, we already noticed that $\widehat{A}$ is a $^\ast$-algebra.  It is easy to check that $\widehat{\Delta} (\omega^\ast) = \widehat{\Delta}(\omega)^\ast$ for all $\omega \in \widehat{A}$ so that $\widehat \Delta$ turns out to be a $^*$-map.
\nl
The next step in showing that the dual $(\widehat A,\widehat \Delta)$, as introduced above, is indeed again an algebraic quantum hypergroup, is the construction of the counit $\widehat\varepsilon$ on $(\widehat A,\widehat \Delta)$. It is just as expected.

\inspr{3.8} Definition \rm 
Let $\omega \in \widehat{A}$ and assume $\omega = \varphi(\,\cdot\, a)$ with $a$ in $A$. Define $\widehat\varepsilon(\omega)=\varphi(a)$. 
\einspr

Also when we use the other expressions for elements in $\widehat A$, we get the expected formulas. So, if $\omega \in \widehat{A}$ is represented as 
$$\omega = \varphi(a \,\cdot\,) = \varphi(\,\cdot\, b) = \psi (c \,\cdot\,) = \psi (\,\cdot\, d)$$
with $a$, $b$, $c$ and $d$ uniquely determined in $A$, then we get 
$$\widehat\varepsilon(\omega) = \varphi(a) = \varphi (b) = \psi (c) = \psi (d).$$
To show that this is correct, one can e.g.\ use the fact that there exists an element $e \in A$ such that $ae = a$ and $eb = b$ (see Proposition 1.6).  Thus, we have $\varphi(a) = \varphi(ae) = \varphi(eb) = \varphi(b)$.
\snl

Then $\widehat\varepsilon$ is a counit on $(\widehat A,\widehat \Delta)$ as follows from the following proposition.

\inspr{3.9} Proposition \rm 
We have that $\widehat{\varepsilon} : \widehat{A} \rightarrow \Bbb C$ is an algebra homomorphism satisfying
\snl
(1) $(\iota \otimes \widehat{\varepsilon}) ((\omega_1 \otimes 1) \widehat{\Delta} (\omega_2)) = \omega_1 \omega_2$ \newline
(2) $(\widehat{\varepsilon} \otimes \iota) (\widehat{\Delta} (\omega_1)(1 \otimes \omega_2)) = \omega_1 \omega_2$
\snl
for all $\omega_1, \omega_2 \in \widehat{A}$.
\snl\bf Proof: \rm  
To prove that $\widehat{\varepsilon}$ is an algebra homomorphism, we put $\omega_1 = \varphi(a \,\cdot\,)$ and $\omega_2 = \varphi(b\,\cdot\,)$.  Then we have $\omega_1 \omega_2 = \varphi(c\,\cdot\,)$ with $c = (\varphi \otimes \iota) (S  \otimes \iota) (\Delta(b) (S^{-1} (a) \otimes 1))$ (see formula (2) in Proposition 3.3).  Therefore, if $\psi = \varphi \circ S$ we have
$$\align
\widehat{\varepsilon} (\omega_1 \omega_2) = \varphi(c) &= \psi ((\iota \otimes \varphi)(\Delta(b) (S^{-1} (a) \otimes 1)))\\
&= \varphi(b) \varphi(a) = \widehat{\varepsilon} (\omega_1) \widehat{\varepsilon} (\omega_2).
\endalign$$
To prove the formula (1), we write $\omega_1 = \psi(a\,\cdot\,)$ and $\omega_2 = \psi (b\,\cdot\,)$.  Then we have $(\omega_1 \otimes 1) \widehat{\Delta}(\omega_2) = (\psi \otimes  \psi) ((\iota \otimes S^{-1}) (\Delta(a) (1 \otimes S(b))) \,\cdot\,)$ (cf.\ the proof of Lemma 3.6). Therefore, we obtain (using formula (4) in Proposition 3.3)
$$
(\iota \otimes \widehat{\varepsilon})((\omega_1 \otimes 1)\widehat{\Delta} (\omega_2)) = \psi ((\iota \otimes \psi) (\iota \otimes S^{-1}) (\Delta(a) (1 \otimes S(b)))\,\cdot\,)
= \omega_1 \omega_2.$$
The formula in (2) is proven in a similar way, now  considering $\omega_1 = \varphi(\,\cdot\, a)$ and $\omega_2 = \varphi(\,\cdot\, b)$ (again, see the proof of Lemma 3.6).\hfill $\blacksquare$
\einspr

We are almost ready to show that $(\widehat{A}, \widehat{\Delta})$ is an algebraic quantum hypergroup in the sense of Definition 1.10.  We first need to define a left integral on $\widehat{A}$.

\iinspr{3.10} Definition \rm 
Let $\psi$ be a right integral on $A$.  For $\omega = \psi(a\,\cdot\,)$, we set $\widehat{\varphi} (\omega) = \varepsilon(a)$.
\einspr

Observe that we use the right integral $\psi$ on $A$ to define the left integral $\widehat\varphi$ on $\widehat A$. In [VD2], it was done the other way. However, as we have used the left integral in the definition of an algebraic quantum hypergroup, we need to define $\widehat\varphi$ first.
\snl  
So, here is the main result of this section ({\it duality for algebraic quantum hypergroups}).

\iinspr{3.11} Theorem \rm Let $(A,\Delta)$ be an algebraic quantum hypergroup. Let the dual 
$(\widehat{A},\widehat{\Delta})$ be defined as before in this section. Then $(\widehat{A},\widehat{\Delta})$ is again an algebraic quantum hypergroup. Moreover, if $(A,\Delta)$ is a $^\ast$-algebraic quantum hypergroup, then $(\widehat{A}, \widehat{\Delta})$ is also a $^\ast$-algebraic quantum hypergroup.

\snl\bf Proof: \rm 
We have already shown that $\widehat A$ is an algebra with a non-degenerate product and that $\widehat \Delta$ is a regular coproduct on $\widehat A$. We also have obtained a counit $\widehat\varepsilon$. 
\snl
We have defined $\widehat\varphi$ above and now, we need to show that this is a left integral. It is clearly non-zero. To show that $\widehat{\varphi}$ is left-invariant on $\widehat{A}$, take $\omega_1$ and $\omega_2$ in $\widehat{A}$ and  calculate $(\iota \otimes \widehat{\varphi}) ((\omega_1 \otimes 1) \widehat{\Delta} (\omega_2))$. Assume $\omega_1 = \psi (a\,\cdot\,)$ and $\omega_2 = \psi (b\,\cdot\,)$ with $a,b\in A$. As in the proof of Lemma 3.6, we get
$$
\langle(\omega_1 \otimes 1)\widehat{\Delta}(\omega_2), x \otimes y \rangle = (\psi \otimes \psi)((\iota \otimes S^{-1})(\Delta(a) (1 \otimes S(b))) (x \otimes y))
$$
for all $x,y \in A$. Therefore, we obtain
$$\align 
(\iota \otimes \widehat{\varphi})((\omega_1 \otimes 1) \widehat{\Delta}(\omega_2)) &= \psi ((\iota \otimes \varepsilon)(\iota \otimes S^{-1}) (\Delta(a) (1 \otimes S(b)))\,\cdot\,)\\
&= \varepsilon(b) \psi (a \,\cdot\,) = \widehat{\varphi} (\omega_2) \omega_1.
\endalign$$
Observe that we used that $\varepsilon$ is invariant under the antipode (see Proposition 2.1).
From the calculation above, we see that $(\iota \otimes \widehat{\varphi}) \widehat{\Delta}(\omega_2)$ and $\widehat{\varphi}(\omega_2)1$ are equal as right multipliers.  So, they also are equal as (two-sided) multipliers in $M(\widehat{A})$. This proves that $\widehat{\varphi}$ is a left integral on $\widehat{A}$.
\snl
Next, we prove that $\widehat{\varphi}$ is faithful.  If $\omega_1$ and $\omega_2$ are elements in $\widehat A$ and if we assume $\omega_1 =  \psi (a\,\cdot\,)$ with $a\in A$, we have
$\omega_1 \omega_2 = \psi ((\iota \otimes \omega_2)(\iota \otimes S^{-1}) \Delta(a)\,\cdot\,)$ as in formula (4) of Proposition 3.3.  Therefore, $\widehat{\varphi}(\omega_1\omega_2) = \omega_2 (S^{-1}(a))$. If this is $0$ for all $a$, then $\omega_2=0$, while if this is $0$ for all $\omega_2$ then $a=0$ (because $\widehat A$ separates points of $A$). This proves the faithfulness of $\widehat\varphi$.  
\snl
Finally, we show that there is an antipode relative to $\widehat{\varphi}$.  For all $\omega \in \widehat{A}$, we define $\widehat{S}(\omega) = \omega \circ S$ on $A$.  It is easy to see that $\widehat{S}(\omega) \in \widehat{A}$. To prove that $\widehat{S}$ is an anti-isomorphism on $\widehat{A}$,  one uses that the antipode flips the coproduct (again see Proposition 2.1). Then, it remains to show that $\widehat{\varphi}$ satisfies the antipode property 
$$\widehat{S} (\iota \otimes \widehat{\varphi})(\widehat{\Delta}(\omega_1) (1 \otimes \omega_2)) = (\iota \otimes \widehat{\varphi}) ((1 \otimes \omega_1) \widehat{\Delta}(\omega_2))$$
for all $\omega_1, \omega_2 \in \widehat{A}$ as in Definition 1.9.
\snl
To prove this, again write $\omega_1 = \psi (a\,\cdot\,)$ and $\omega_2 = \psi (\,\cdot\, b)$ with $a,b\in A$.  Then we have $(1 \otimes \omega_1) \widehat{\Delta}(\omega_2) = \sum \limits_i \psi (\,\cdot\, q_i) \otimes \psi (p_i \,\cdot\,)$ where $\sum\limits_i p_i \otimes q_i = (\iota \otimes S^{-1})((1 \otimes S(b)) \Delta(a))$ as in the proof of Proposition 3.7. The right hand side of the antipode equation above is
$$(\iota \otimes \widehat{\varphi}) ((1 \otimes \omega_1) \widehat{\Delta}(\omega_2))=\sum\limits_i \varepsilon(p_i) \psi(\,\cdot\, q_i) = \psi (\,\cdot\, S^{-1} (a) b).$$
For the left hand side, we first calculate the expression for $\widehat{\Delta} (\omega_1) (1 \otimes \omega_2)$ in $\widehat{A} \otimes \widehat{A}$, using the given representations for $\omega_1$ and $\omega_2$.  For all $x$, $y$ in $A$, we have
$$\align
\langle \widehat{\Delta} (\omega_1) (1 \otimes \omega_2), x \otimes y\rangle &= \langle\omega_1 \otimes \omega_2, (x \otimes 1) \Delta(y)\rangle\\
&= \omega_2 ((\psi \otimes \iota) (ax \otimes 1) \Delta (y)) \\
&= \psi ((\iota \otimes (\omega_2 \circ S^{-1})) (\Delta (ax) (y \otimes 1)))\\
&= \psi ((\iota \otimes (\omega_2 \circ S^{-1})) (\Delta(ax)) y).
\endalign$$
Therefore, $\langle(\iota \otimes \widehat{\varphi}) (\widehat{\Delta} (\omega_1) (1 \otimes  \omega_2)), x \rangle = \varepsilon (\iota \otimes (\omega_2 \circ S^{-1})) (\Delta(ax)) = (\omega_2 \circ S^{-1}) (ax)$.
So, we have $(\iota \otimes \widehat{\varphi})(\widehat{\Delta}(\omega_1)(1\otimes \omega_2)) = (\omega_2 \circ S^{-1}) (a \,\cdot\,)$ and $\widehat{S} (\iota \otimes \widehat{\varphi}) (\widehat{\Delta} (\omega_1) (1 \otimes \omega_2)) = \omega_2 (\,\cdot\, S^{-1} (a)) = \psi(\,\cdot\, S^{-1}(a)b)$. We see that the left and the right hand side of antipode equation above are equal. 
\snl
This completes the proof of the fact that $(\widehat A,\widehat \Delta)$ is an algebraic quantum hypergroup. 
\snl
Finally, if $(A,\Delta)$ is a $^\ast$-algebraic quantum hypergroup, we already mentioned that $\widehat{A}$ is a $^\ast$-algebra and $\widehat{\Delta}$ is a $^\ast$-map. It follows that also $(\widehat A, \widehat \Delta)$ is a $^*$-algebraic quantum hypergroup. \hfill$\blacksquare$
\einspr

If we set $\widehat\psi=\widehat\varphi\circ \widehat S$ as we do for the original pair $(A,\Delta)$ we find easily that $\widehat\psi(\omega)=\varepsilon(a)$ when $\omega=\varphi(\,\cdot\, a)$. We will use this formula in the proof of Theorem 3.12 below.
\snl
In the case of a $^*$-algebraic quantum hypergroup, as we mentioned earlier (see a remark following Definition 1.5), it makes sense to assume that the integrals are positive. Now, suppose that $\psi$ is positive on $A$. Then, it can be shown that $\widehat\varphi$ is again positive on $\widehat A$. Indeed, assume that $\omega=\psi(a\,\cdot\,)$ with $a\in A$. Then, as we see from the fourth formula in Proposition 3.3, we get $\omega\omega^*=\psi(e\,\cdot\,)$ where $e=(\iota\ot(\omega^*\circ S^{-1}))\Delta(a)$. Therefore
$$\widehat\varphi(\omega\omega^*)=\varepsilon(e)=\omega^*(S^{-1}(a))=\omega(a^*)^-=\psi(aa^*)^-.$$
So, we see that indeed $\widehat\varphi$ is positive when $\psi$ is positive. Remark that there seems to be no obvious way to show that a positve left integral exists when there is a positive right integral. The result is known to be true for $^*$-algebraic quantum groups (see [K-VD]), but the proof is quite involved. It has to be investigated if the result is still true for algebraic quantum hypergroups. We refer to Section 5 for a further discussion about this problem.
\nl
We will illustrate the main result at the end of this section using the motivating example from Section 1 (Example 1.11). And in the next section, we will consider the objects associated with the dual $(\widehat A,\widehat\Delta)$ (the modular element $\widehat\delta$, the modular automorphisms $\widehat\sigma$ and $\widehat\sigma'$ and the scaling constant of the dual) and see how they can be found from the data of the original algebraic quantum hypergroup $(A,\Delta)$.
\nl
Now, we will show that taking the dual of $(\widehat A,\widehat\Delta)$ will give us back the original pair $(A,\Delta)$. This is the content of the following theorem ({\it biduality for algebraic quantum hypergroups}).

\iinspr{3.12} Theorem \rm
Let $(A, \Delta)$ be an algebraic quantum hypergroup.  Let $(\widehat{A}, \widehat{\Delta})$ be the dual algebraic quantum hypergroup. For $a \in A$ and $\omega \in \widehat{A}$, we set $\Gamma(a)(\omega) = \omega(a)$.  Then $\Gamma(a) \in \dubbeldualA$ for all $a\in A$. Moreover, $\Gamma$ is an isomorphism between the algebraic quantum hypergroups $(A,\Delta)$ and $(\dubbeldualA, \widehat{\widehat{\Delta}})$.  In the $^\ast$-case, we have that $\Gamma$ is a $^\ast$-isomorphism.

\snl\bf Proof: \rm  For $a$ in $A$, define $\Gamma(a)$ as a linear functional on $\widehat A$ by $\langle\omega,\Gamma(a)\rangle = \omega(a)$ whenever $\omega\in \widehat A$. We will first show that actually, $\Gamma(a)\in \dubbeldualA$. In order to prove this, denote $\omega=\varphi(\,\cdot\, S(a))$ and take any $\omega_1$ in $\widehat{A}$. Then we get, by using formula (1) of Proposition 3.3, that $\omega_1\omega = \varphi(\,\cdot\, d)$ where
$$d = ((\omega_1 \circ S^{-1}) \otimes  \iota) \Delta(S(a)) = S ((\iota\otimes \omega_1) \Delta(a)).$$
Therefore $\widehat{\psi} (\omega_1\omega) = \varepsilon(d) = \omega_1(a) = \langle\omega_1,\Gamma(a)\rangle$ and thus we have $\Gamma(a) = \widehat{\psi} (\,\cdot\, \omega)$ and $\Gamma(a)\in \dubbeldualA$. 
\snl
It is clear that $\Gamma$ is an isomorphism between the linear spaces $A$ and $\dubbeldualA$. That $\Gamma$ respects the multiplication and the comultiplication is straightforward because in both cases, the product is dual to the coproduct and vice versa. So $\Gamma$ is an isomorphism of algebraic quantum hypergroups.
\snl
In the case of a $^*$-algebra, it is easily proven that $\Gamma$ is also a $^*$-isomorphism.
\hfill $\blacksquare$
\einspr

By uniqueness of the counit, the antipode and the integrals of an algebraic quantum group, one must have that $\Gamma$ also respects these objects. For the counits e.g.\, we get
 $$\dubbeldualepsilon (\Gamma(a)) = \dubbeldualepsilon (\widehat{\psi} (\,\cdot\,  \omega)) =  \widehat{\psi} (\omega) = \varepsilon(S(a)) = \varepsilon(a)$$
where, as before, $a\in A$ and $\omega=\varphi(\,\cdot\,S(a))$. For the antipodes, we get the result essentially because the antipode on $\widehat A$ is defined as the adjoint of the antipode on $A$ (see the proof of Theorem 3.11). 
Finally, for the integrals, it is somewhat more complicated. One can show that also here 
$\dubbeldualphi (\Gamma(a))=\varphi(a)$ for all $a$ when we use the conventions $\psi=\varphi\circ S$ and $\widehat\psi=\widehat \varphi\circ \widehat S$ and when $\widehat\varphi$ is defined using $\psi$ as in Definition 3.10. To prove this result however, we need a formula for the modular automorphism of $\widehat\psi$. This formula is proven in the next section and there we will also give the argument for the equality $\dubbeldualphi (\Gamma(a))=\varphi(a)$ (see a remark following Proposition 4.1). 
\nl
Let us again finish this section by looking at some {\it special cases} and {\it examples}.
\snl
First consider the case of an algebraic quantum hypergroup $(A,\Delta)$ of compact type (cf.\ Definition 1.13). So, $A$ has an identity $1$ and $\Delta(1)=1 \otimes 1$. Then $\varphi\in \widehat A$ and from formula (1) in Proposition 3.3, we easily find that $\omega\varphi=\widehat\varepsilon(\omega)\varphi$ for all $\omega\in \widehat A$. Similarly $\psi\in \widehat A$ and $\psi\omega=\widehat\varepsilon(\omega)\psi$ for all $\omega\in \widehat A$. This will follow from formula (3) in the same proposition.
\snl
This takes us to the following notion.

\iinspr{3.13} Definition \rm
Let $(A,\Delta)$ be an algebraic quantum hypergroup with counit $\varepsilon$. An element $h \in A$ is called a left co-integral if it is non-zero and if $ah=\varepsilon(a)h$ for all $a\in A$. Similarly, a right co-integral is a non-zero element $k\in A$ so that $ka=\varepsilon(a)k$ for all $a\in A$.
\einspr

It is clear that a left co-integral exists if an only if a right co-integral exists (apply the antipode). Then we come to the following definition.

\iinspr{3.14} Definition \rm
An algebraic quantum hypergroup is called of {\it discrete type} if there exists a left co-integral.
\einspr

We again refer to the discussions in the last section to clarify this terminology.
\snl
In the remark, preceding the Definition 3.13, we saw that $(\widehat A,\widehat \Delta)$ has co-integrals when $A$ has an identity. In fact, we have the following result.

\iinspr{3.15} Proposition \rm
Let $(A,\Delta)$ be an algebraic quantum group and $(\widehat A,\widehat \Delta)$ its dual. Then $(A,\Delta)$ is of compact type if and only if $(\widehat A,\widehat\Delta)$ is of discrete type.
\snl\bf Proof: \rm
We have already given an argument for one direction. We will now prove the converse, but for the dual. The result will then follow from (bi)duality.
\snl
So assume that $(A,\Delta)$ is of discrete type and let $h$  be a left co-integral. For all $a$ in $A$ we have $\varphi(ah)=\varepsilon(a)\varphi(h)$. Because $\varphi$ is faitful and $h\neq 0$, we must have $\varphi(h)\neq 0$. Then $\varepsilon\in\widehat A$ and so $\widehat A$ is unital. Therefore, it is of compact type.\hfill $\blacksquare$
\einspr
  
It is not so hard to show that, if $(A,\Delta)$ is both of discrete and of compact type, it must be finite-dimensional. Indeed, assume that $h$ is a left co-integral. We have seen that $\varphi(h)\neq 0$ and so we can assume that $\varphi(h)=1$. Then, for all $a$ in $A$ we have $\Delta(a)(1\ot h)=(\iota\ot\varepsilon)\Delta(a) \ot h= a\ot h$ and so
$(\iota\ot\varphi)(\Delta(a)(1\ot h))=a$. By the antipode property we get $S(a)=(\iota\ot\varphi)((1\ot a)\Delta(h))$ and we see that $A$ is part of the 'left leg' of $\Delta(h)$. If however $1\in A$, this is a finite-dimensional subspace of $A$ and therefore $A$ itself must be finite-dimensional. 
\snl
Also conversely, when $A$ is finite-dimensional, it must be of compact and of discrete type. One possible argument is simple. By Proposition 1.6 we know that $A$ has local units. Because $A$ is finite-dimensional, it must have a unit. Similarly for the dual $\widehat A$. 
\snl
When $A$ is finite-dimensional, we call $(A,\Delta)$ of {\it finite type}.
\nl
We finally look again at Example 1.11.

\iinspr{3.16} Example \rm
Remember that $G$ is a group and $H$ a finite subgroup. The algebra $A$ is the space of complex functions with finite support on $G$ and constant on double $H$-cosets. The product is pointwise and the coproduct $\Delta$ is defined by
$$\Delta(f)(p,q)=\frac1n \sum_{h\in H} f(phq)$$
where $n$ is the number of elements in $H$ and $p,q\in G$ and $f\in A$. The counit $\varepsilon$ is given by $\varepsilon(f)=f(e)$ where $e$ is the identity of $G$. The antipode $S$ is given by $S(f)(p)=f(p^{-1})$ for $f\in A$ and $p\in G$. The left integral $\varphi$ is given by 
$$\varphi(f)=\sum_{p\in G}f(p)$$
and the right integral $\psi$ is equal to $\varphi$.
\snl
It follows immediately from the definition of the dual that $\widehat A$ is realized also as the space of complex funtions on $G$ with finite support and constant on double $H$-cosets with the pairing 
$$\langle f,g\rangle = \sum_{p\in G} f(p)g(p)$$
for $f\in A$ and $g\in \widehat A$. The product in $\widehat A$, dual to the coproduct on $A$ is easily calculated and we get the (ordinary) convolution product
$$(g_1g_2)(p)=\sum_{q\in G} g_1(q)g_2(q^{-1}p)$$
when $g_1,g_2\in \widehat A$ and $p\in G$. On the other hand, for the coproduct $\widehat \Delta$ on $\widehat A$, we get $\widehat\Delta(\pi(p))=\pi(p)\ot \pi(p)$ where, for $p\in G$, we let
$$\pi(p)=\frac{1}{n^2}\sum_{h,h'\in H} \lambda_{hph'}$$
and where $\lambda_q$ is the function on $G$ that is one in $q$ and $0$ everywhere else. The counit $\widehat\varepsilon$ on $\widehat A$ is given by $\widehat\varepsilon(\pi(p))=1$ for all $p$. The antipode $\widehat S$ on $\widehat A$ is given by $\widehat S(\pi(p))=\pi(p^{-1})$. The dual left integral $\widehat\varphi$ is given by $\widehat\varphi(\pi(p))=0$ except when $\pi(p)=\pi(e)$. The right integral $\widehat\psi$ is again equal to the left integral $\widehat\varphi$.
\einspr

In fact, there is a better way to look at this example. Indeed, let $B$ be the algebra of all complex functions on $G$ with finite support, with the convolution product and with the coproduct given by $\Delta(\lambda_p)=\lambda_p \otimes \lambda_p$. Consider the element $u$ in $B$, given by $u=\frac{1}{n}\sum_{h\in H} \lambda_h$ (which is nothing else but the element $\pi(e)$ as above). Then $u^2=u$ and also 
$$\Delta(u)(1\ot u)=u\ot u.$$
This means that $u$ is a so-called group-like projection (in the sense of Definition 1.1 in [L-VD1]). Then $\widehat A=uBu$ and $\widehat\Delta(b)=(u\ot u)\Delta(b)(u\ot u)$ for any $b\in \widehat A$. The counit, the antipode and the integrals are simply the restrictions to $\widehat A$ of resp.\ the counit, the antipode and the integrals on $B$. We get a typical example of an algebraic quantum hypergroup of compact type as in Section 2 of [L-VD1]. The identity in $\widehat A$ is noting else but $u$.
\snl
Of course, Example 1.11 is a typical example of an algebraic quantum hypergroup of discrete type. A left co-integral is the function that is $1$ on $H$ and $0$ everywhere else. In this case, a left co-integral is also a right co-integral.
\snl
For less trivial examples, we again refer to [D-VD2]. 
\nl\nl

\bf 4. More properties of an algebraic quantum hypergroup and its dual \rm
\nl
In this section, we will collect some more results and formulas, involving the data associated with an algebraic quantum hypergroup $(A,\Delta)$ as well as of its dual $(\widehat A,\widehat\Delta)$. We also consider some module structures involving both $A$ and its dual $\widehat A$.
\nl
So, in what follows, we have an algebraic quantum hypergroup $(A,\Delta)$, with counit $\varepsilon$, antipode $S$ and left and right integrals $\varphi$ and $\psi$ respectively. We assume that $\psi=\varphi\circ S$. The modular element relating the left and the right integral is $\delta$ and the modular automorphisms associated with $\varphi$ and $\psi$ respectively  are $\sigma$ and $\sigma'$. Finally, there is the scaling constant $\tau$. We have all the relations among these objects as obtained in Section 2. Now, because the dual $(\widehat A,\widehat\Delta)$ is again an algebraic quantum hypergroup, we also have the data associated with this dual. The objects here are denoted as for $(A,\Delta)$ but with a {\it hat}. We normalize $\widehat \varphi$ by using the Definition 3.10. So $\widehat\varphi(\omega)=\varepsilon(a)$ when $\omega=\psi(a\,\cdot\,)$ with $a\in A$. Again, we normalize $\widehat\psi$ by $\widehat\psi=\widehat\varphi\circ\widehat S$. Then, as we have seen before, we get $\widehat \psi(\omega)=\varepsilon(a)$ when $\omega=\varphi(\,\cdot\, a)$ with $a$ in $A$.
\snl
We have the following formulas for the dual objects, in terms of the data associated with the original quantum hypergroup. Observe that the formulas are the same as in the case of an algebraic quantum group (see [Ku] and e.g.\ also [D-VD1]).

\inspr{4.1} Proposition \rm 
The modular element $\widehat \delta$ and its inverse ${\widehat\delta}^{-1}$, when considered as linear functionals on $A$, are given by 
$$\align  &\widehat{\delta} = \varepsilon \circ \sigma^{-1} = \varepsilon \circ {\sigma'}^{-1} \\
 &\widehat{\delta}^{-1} = \varepsilon \circ \sigma = \varepsilon \circ \sigma'.
\endalign$$
On the other hand, the modular automorphisms $\widehat\sigma$ and $\widehat\sigma'$, associated with $\widehat\varphi$ and $\widehat\psi$ respectively, satisfy
$$\align \langle\widehat{\sigma} (\omega),a \rangle &= \langle\omega, S^2(a)\delta^{-1}\rangle \\
         \langle\widehat{\sigma}' (\omega), a\rangle &= \langle\omega, \delta^{-1}S^{-2}(a)\rangle
\endalign$$
for all $a\in A$ and $\omega\in\widehat A$.

\snl\bf Proof: \rm 
We first prove the formulas for the modular element $\widehat\delta$. Take $\omega_1, \omega_2$ in $\widehat{A}$. Using Proposition 2.5, we see that we must have 
$$\text{(1)}\qquad\qquad(\widehat{\varphi} \otimes \iota)(\widehat{\Delta} (\omega_1)(1\otimes \omega_2)) = \widehat{\varphi} (\omega_1) \widehat{\delta} \omega_2.$$
To calculate the left hand side of this equation (1), take $\omega_1 = \psi(a\,\cdot\,)$ and $\omega_2 = \varphi(\,\cdot\, b)$ with $a,b\in A$. Then we have, for all $x, y \in A$, that
$$\align
\langle \widehat{\Delta} (\omega_1)(1 \otimes \omega_2), x \otimes y\rangle &= \langle \omega_1 \otimes \omega_2, (x \otimes 1) \Delta(y)\rangle\\
&= \psi (ax S^{-1} (\iota\otimes \varphi)((1 \otimes y)\Delta(b)))\\
&= \psi(({\sigma'}^{-1} \circ S^{-1}) (\iota \otimes \varphi) ((1 \otimes y) \Delta(b)) ax).
\endalign$$
Therefore we get,
$$\align
\langle (\widehat{\varphi} \otimes \iota) (\widehat{\Delta} (\omega_1)(1 \otimes \omega_2)), y\rangle
&= \varepsilon(a)(\varepsilon \circ {\sigma'}^{-1} \circ S^{-1})((\iota \otimes \varphi) ((1 \otimes y)\Delta(b))) \\
&= \widehat{\varphi} (\omega_1) \varphi(y((\varepsilon \circ {\sigma'}^{-1} \circ S^{-1}) \otimes \iota) \Delta(b)).
\endalign$$
Next, we calculate the right hand side of equation (1).  Considering $\widehat\delta$ as a linear functional on $A$ and using Proposition 3.4, we have that, for all $y \in A$,
$$\widehat{\varphi}(\omega_1) \langle \widehat{\delta}\omega_2, y\rangle = \widehat{\varphi} (\omega_1) \varphi(y(((\widehat{\delta}  \circ S^{-1}) \otimes \iota) \Delta(b))).$$
Comparing these two expressions for the left and the right hand side of equation (1) for all $\omega_1$ and $\omega_2$ in $\widehat{A}$ and using that the left integral $\varphi$ on $A$ is faithful, we conclude that, for all $b\in A$ we have
$$
((\varepsilon \circ {\sigma'}^{-1} \circ S^{-1}) \otimes \iota)  \Delta(b) = ((\widehat{\delta} \circ S^{-1}) \otimes \iota) \Delta(b).$$
If now, we apply the counit $\varepsilon$ to both sides of this equation, we  obtain
$$ (\varepsilon \circ {\sigma'}^{-1} \circ S^{-1}) (b) = (\widehat{\delta} \circ S^{-1})(b)$$
for all $b$ in $A$. Hence, we see that $\widehat{\delta} = \varepsilon \circ {\sigma'}^{-1}$ as linear functionals on $A$.
\snl
Recall from Proposition 2.7 that  $\delta \sigma (a) = \sigma' (a) \delta$ for all $a \in A$.  Also
$\delta {\sigma'}^{-1} (a) = \sigma^{-1} (a) \delta$ for all $a \in A$. Because $\varepsilon(\delta)=1$, we get 
$\varepsilon \circ \sigma = \varepsilon \circ \sigma'$ as well as $\varepsilon \circ\sigma^{-1}  = \varepsilon \circ {\sigma'}^{-1}$. On the other hand, using the other relation in Proposition 2.7, namely $\sigma S \sigma'=S$, we find the formulas for the inverse $\widehat\delta^{-1}$ by applying the antipode.
\snl
Let us now prove the formulas for the modular automorphisms $\widehat\sigma$ and ${\widehat\sigma}'$. The automorphism $\widehat{\sigma}$ is characterized by  $\widehat{\varphi} (\omega_1\omega_2) = \widehat{\varphi}({\omega_2} \widehat{\sigma} (\omega_1))$ when $\omega_1, \omega_2 \in \widehat{A}$. Now take $\omega_1 = \psi (a\,\cdot\,)$ and $\omega_2 = \psi(b\,\cdot\,)$ with $a,b\in A$. Using the formula (4) in Proposition 3.3, we have $\omega_1 \omega_2 = \psi(d\,\cdot\,)$ where $d = (\iota \otimes (\omega_2 \circ S^{-1}))\Delta(a)$.  Therefore, 
$$\widehat{\varphi} (\omega_1\omega_2) = \varepsilon(d) = \psi (bS^{-1} (a)) = \psi(S^{-1} (aS(b))) = \psi (aS(b) \delta^{-1}) = \omega_1 (S(b) \delta^{-1}).$$ 
Now assume that $\omega_3 = \psi(c\,\cdot\,)$ with $c\in A$. We have, in a similar way as above, that $\widehat{\varphi}(\omega_2\omega_3) = \psi (c S^{-1}(b))$. If we want to have that $\omega_3=\widehat\sigma(\omega_1)$, we need the element $c$ to satisfy
$$\omega_1(S(b)\delta^{-1}) =\psi (c S^{-1}(b))$$
for all $b\in A$. So, we must have
$\omega_1 (S^2(b \delta^{-1})) = \psi(cb)$ for all $b$ and we see that $\widehat{\sigma} (\omega_1)(b)  = (\omega_1 \circ S^2) (b\delta^{-1})$.  
\snl
The formula for the automorphism $\widehat{\sigma}'$ can be found in a similar way, but it is easier to deduce it from the formula for $\widehat\sigma$, using also that $\widehat\sigma \circ S\circ \widehat\sigma' =S$ . \hfill $\blacksquare$
\einspr

We have the following interesting consequence. Because $\widehat{\delta} = \varepsilon \circ \sigma^{-1}$ on $A$, we see that $\widehat{\delta}$ is an algebra homomorphism on $A$.  Therefore, we have $\langle aa',\widehat\delta \rangle=\langle a,\widehat\delta\rangle \langle a',\widehat\delta  \rangle$. We can interprete this formula as $\widehat{\Delta}(\widehat{\delta}) = \widehat{\delta} \otimes \widehat{\delta}$ in the dual space $(A \otimes A)'$. By duality, we also have the formula $\Delta(\delta)=\delta\ot\delta$ in the dual space $(\widehat A\otimes \widehat A)'$.
\snl
It is more difficult to consider these formulas as valid in $M(\widehat A\otimes \widehat A)$ and $M(A\ot A)$ respectively because for this, we would first need to extend the coproducts and this is more subtle because the coproducts are no longer homomorphisms.
\snl
The formula in Proposition 4.1 can now be used to prove the result, announced in the previous section after the proof of Theorem 3.1 (in connection with biduality). If $\Gamma:A\to \dubbeldualA$ is defined as before an if again $\omega=\varphi(\,\cdot\,S(a))$ we get indeed
$$\align
\dubbeldualphi (\Gamma(a)) &= \dubbeldualphi(\widehat{\psi} (\,\cdot\, \omega)) = \dubbeldualphi (\widehat{\psi} ((\widehat{\sigma}')^{-1} (\omega)\,\cdot\,))
= \widehat{\varepsilon} ((\widehat{\sigma}')^{-1} (\omega))\\
&= \widehat{\varepsilon} ((\omega \circ S^2)(\delta\,\cdot\,)) = \omega(\delta)= \varphi(\delta S(a)) = \varphi S(a \delta^{-1}) = \varphi(a\delta^{-1}\delta)=\varphi(a).
\endalign$$
We will give another interpretation of the second pair of formulas in the previous proposition and relate these formulas to earlier statements after we have properly defined the obvious module structures below.

\nl
\it Module structures \rm
\nl
Now, we will look at the obvious {\it module structures} for an algebraic quantum hypergroup.
\snl
So, as before, let $(A, \Delta)$ and $(\widehat{A}, \widehat{\Delta})$ denote an algebraic quantum hypergroup and its dual.  Just as in the case of ordinary algebraic quantum groups, we consider four module-structures, denoted in the following way
$$A \blacktriangleright \widehat{A} \qquad \widehat{A} \blacktriangleleft A \qquad \widehat{A} \blacktriangleright A \qquad A \blacktriangleleft \widehat{A}.$$ 
Here are the precise definitions.

\inspr{4.2} Definition \rm
For $a \in A$ and $\omega \in \widehat{A}$, these module actions are defined by the formulas
$$\alignat{2}
a \blacktriangleright \omega &= \omega (\,\cdot\, a) &\qquad\qquad \omega \blacktriangleleft a &= \omega(a \,\cdot\,)\\
\omega \blacktriangleright a &= (\iota \otimes \omega) \Delta(a) &\qquad\qquad a \blacktriangleleft \omega &= (\omega \otimes \iota) \Delta(a).
\endalignat$$
\einspr

The formulas in the first row define left and right $A$-modules because the product in $A$ is associative.  The formulas in the second row define left and right $\widehat{A}$-modules because the comultiplication on $A$ is coassociative.  Observe that for all $\omega$, $\omega'$ in $\widehat{A}$ and $a \in A$, we have $\langle \omega', \omega \blacktriangleright a\rangle = \langle \omega' \omega, a\rangle$ and of course also $\langle \omega', a\blacktriangleleft \omega \rangle = \langle \omega\omega',a\rangle$.  If $(A, \Delta)$ is an algebraic quantum group, in the sense that $\Delta$ is also an algebra homomorphism, the module structures are module algebras.
\snl
The modules $A \btr \widehat{A}$ and $\widehat{A} \btl A$ are unital in the sense that $A \btr \widehat{A} = \widehat{A} = \widehat{A} \btl A$.  This follows from the existence of local units as obtained in Proposition 1.6. Similarly, we have $\widehat{A} \btr A = A = A \btl \widehat{A}$ because the integrals on $A$ are faithful and Lemma 1.8 can be applied.  Observe that for $a \in A$, there is an element $\omega \in \widehat{A}$ such that $\omega \btr a = a$.  Indeed, $a = \sum \omega_i \btr a_i$ for $\omega_i \in \widehat{A}$ and $a_i \in A$ and using Proposition 1.6 (applied for the dual), we can choose $\omega \in \widehat{A}$ such that $\omega \omega_i = \omega_i$ for all $i$.  Therefore we have 
$$\omega \btr a = \sum \omega \omega_i \btr a_i = \sum \omega_i \btr a_i = a.$$
From this observation, we see that the module structure $\widehat{A} \btr A$ extends to $M(\widehat{A}) \btr A$ in a natural way.  E.g. $\widehat{\delta} \btr a = \widehat{\delta} \omega \btr a$ where $\omega \in \widehat{A}$ is choosen such that $\omega \btr a = a$.  One easily checks that this definition does not depend on the choice of $\omega$ in $\widehat{A}$.
\snl
We also have extended the pairing between $A$ and $\widehat A$ to a bilinear map on $M(\widehat A) \times A$ (see the remark following Proposition 3.4). The formulas we obtained there can now be rewritten, using the above module actions.
We get
$$\align \langle f \omega , x \rangle &= \langle f, \omega \btr x \rangle \\
         \langle \omega f, x \rangle &= \langle f, x \btl \omega \rangle
\endalign$$
where $\omega\in\widehat A$, $f\in M(\widehat A)$ and $x\in A$. This is completely in accordance with similar formulas considered earlier with also $f\in \widehat A$.
\snl
This observation is also important to get a better understanding of the first two formulas in the next proposition. In this last proposition of this section, we also obtain Radford's formula for the 4th power of the antipode as an easy consequence of the other formulas.
 
\inspr{4.3} Proposition \rm 
For all $a \in A$, we have
$$\sigma(a) = \widehat{\delta}^{-1} \btr S^2(a) \qquad \qquad \sigma' (a) = S^{-2}(a) \btl \widehat{\delta}^{-1}$$
$$S^4(a) = \delta^{-1}(\widehat{\delta}  \btr a \btl \widehat{\delta}^{-1})\delta.$$

\snl\bf Proof: \rm 
Essentially, the first two formulas were encounterd already in the proof of the last statement of Proposition 2.7. Let us consider the argument again for the formula for $\sigma$. The formula for $\sigma'$ can be shown in a similar way, or by using the fact that $\sigma S \sigma'=S$ (see (1) in Proposition 2.7). 
\snl
Now, from (5) in the same proposition, we have $\Delta(\sigma^{-1}(a)) = (S^{-2}\otimes \sigma^{-1}) \Delta(a)$ for all $a \in A$.  In Proposition 4.1, we have proven that $\widehat{\delta} = \varepsilon \circ \sigma^{-1}$ as linear maps on $A$.  Combining these results, we have $\widehat{\delta} \btr a = S^2 (\sigma^{-1}(a))$ for all $a \in A$.  Now the formula of $\sigma(a)$ easily follows.   
\snl
Next, recall from Proposition 2.7 (1) that $\delta \sigma(a) = \sigma' (a) \delta$ for all $a \in A$.  Substituting the above formulas in this equation, leads to the equality
$$\delta(\widehat{\delta}^{-1} \btr S^4(a)) = (a \btl \widehat{\delta}^{-1})\delta$$
for all $a \in A$. This means $S^4(a) = \widehat{\delta} \btr (\delta^{-1} (a \btl \widehat{\delta}^{-1} ) \delta)$.
\snl
Now, because $\sigma^{-1}(\delta)=\tau\delta$ (see Proposition 2.7), we will get $\widehat \delta \btr (x\delta)=S^2\sigma^{-1}(x\delta)=\tau (S^2\sigma^{-1}(x))\delta=\tau (\widehat\delta\btr x)\delta$ for all $x$. Similarly $\widehat \delta \btr(\delta^{-1}x)=\tau^{-1}\delta^{-1}(\widehat\delta\btr x)$ for all $x$. Therefore we have $S^4(a) = \delta^{-1}(\widehat{\delta} \btr a \btl \widehat{\delta}^{-1})\delta$ for all $a$ and this completes the proof. \hfill$\blacksquare$
\einspr
Observe that the first two formulas in this proposition are essentially the same as the last two formulas in Proposition 4.1 (in dual form). The proof here is a bit different and more adapted to this new formulation.
\snl
In [D-VD1, Theorem 1.6], we have obtained the same formulas for $\sigma$, $\sigma'$ and $S^4$ in the case of  an algebraic quantum group. The new results here are more general than the ones obtained earlier in [D-VD1] and the proofs are not more complicated. Radford's formula for $S^4$ is well-known in the case of finite-dimensional Hopf algebras. We see that it can be extended, not only to the case of algebraic quantum groups, but even to algebraic quantum hypergroups. As we already mentioned, it turns out to be an easy consequence of the different relations between the data, as obtained earlier in this paper. See also [D-VD-W] for a discussion about this approach to Radford's formula.
\nl\nl

\bf 5. Conclusions and further research \rm
\nl
In this paper, we have developed the theory of algebraic quantum hypergroups. The definition given in Section 1 (Definition 1.10) is perhaps not the final one. The antipode property (cf.\ Definition 1.9) characterizes the antipode in terms of the left integral. Therefore, it seems (at present) not possible to extend the notion of a (multiplier) Hopf algebra first to the case where the coproduct is no longer assumed to be an algebra homomorphism and later define algebraic quantum hypergroups as such objects carrying integrals. This is how algebraic quantum groups were introduced in [VD1] and [VD2].
\snl
On the other hand, we are convinced that we have the right concept. The fact that duality works fine is a strong indication for this statement. Moreover, and as we mentioned already earlier, this is quite remarkable, all of the data and most of the relations among these data that are known for ordinary algebraic quantum groups, remain present for these hypergroups. Furthermore, it turned out that the proofs in this more general setting were not more difficult than for algebraic quantum hypergroups. On the contrary, some of the arguments are even simpler!
\snl
We have discussed a few examples. The basic one is Example 1.11, constructed from a finite subgroup $H$, not necessarily a normal subgroup, of a group $G$. There is also the dual of this example (cf.\ Example 3.16). Both are special cases of a more general situation that we encounter in a paper on discrete and compact subgroups of algebraic quantum groups (cf.\ [L-VD1]). In fact, that paper is what motivated us to start the study of these algebraic quantum hypergoups as these objects arose naturally in [L-VD1]. One of the examples is of discrete type (Example 1.11) while the other is of compact type (Example 3.16). Unfortunately, these examples are far too simple to illustrate the many nice features of the algebraic quantum hypergroups as developed in this paper. On the other hand, this paper is already quite long and therefore, we have chosen to give the more complicated examples in a separate paper (cf.\ D-VD2]). They will enable us to illustrate the various data and their relations as obtained in this paper.
\snl
In the case of a $^*$-algebraic quantum group with positive integrals, it is possible to represent the underlying $^*$-algebra as a $^*$-algebra of bounded operators on a Hilbert space. Doing so in [K-VD], it is shown that any such $^*$-algebraic quantum group can be 'completed' to a C$^*$-algebraic quantum group (or in more modern language, to a locally compact quantum group as in [K-V1] and [K-V2]). In fact, the development of the theory of locally compact quantum groups was greatly inspired by the theory of $^*$-algebraic quantum groups (as developed in [VD1] and [VD2]) and the passage to a C$^*$-algebra (as in [K-VD]). Therefore, it is an obvious question whether the same can be done for $^*$-algebraic quantum hypergroups with positive integrals. Is it possible also here to represent the underlying $^*$-algebra as a $^*$-algebra of bounded operators on a Hilbert space and can the coproduct be extended by continuity to the closure of this algebra? If this is possible, this step should eventually yield the development of what could be called a locally compact quantum hypergroup.
\snl
Related with this is the following, non-trivial problem. If we have a positive left integral $\varphi$ on a $^*$-algebraic quantum hypergroup $(A,\Delta)$, is there also a positive right integral? This is known to be the case for $^*$-algebraic quantum groups. It was shown in [K-VD], but the result is highly non-trivial. In general, there seems to be no reason why the natural candidate for the right integral, namely $\varphi\circ S$, should be positive when $\varphi$  is positive. The fact that $S$ is not a $^*$-map causes the problem. Roughly speaking, what is needed is a 'self-adjoint' square root $\delta^\frac12$ of the modular element $\delta$ in the multiplier algebra $M(A)$. Then one could put $\psi=\varphi(\delta^\frac12\,\cdot\, \delta^\frac12)$. Another possible solution is obtained from a so-called polar decomposition of the antipode. One of the ingredients is a $^*$-anti-automorphism that flips the coproduct. This will convert a positive left integral into a positive right integral. See a recent paper [DC-VD] where a new treatment of this problem is considered.
We strongly believe that the result is true also for algebraic quantum hypergroups. The fact that a positive left integral on the dual $(\widehat A,\widehat \Delta)$ exists when there is a positive right integral on $(A,\Delta)$ is a strong indication (as well as the fact that the result holds for algebraic quantum groups).
\snl
There is also an existing theory of compact quantum hypergroups (see [C-V]). It is expected that the above procedure, if applied to a $^*$-algebraic quantum hypergroup of compact type with positive integrals, will yield a compact quantum hypergroup in the sense of [C-V].
\snl
Also the terminology introduced in this paper ('compact type' and 'discrete type') refers to the point of view that the algebraic quantum hypergroups are seen as a step towards a theory of locally compact quantum hypergroups. They should be, in the first place, quantizations of locally compact spaces in the sense that the underlying algebras should be non-commutative analogues of the C$^*$-algebra of continuous complex functions tending to $0$ at infinity on these spaces. The underlying quantum space is said to be compact if the algebra has an identity. This is where the notion 'compact type' comes from in this paper (cf.\ Definition 1.13). The use of 'discrete type' can be motivated in a similar way, but it is best understood if we refer to the duality of Pontryagin for locally compact abelian groups. There, it is known that the dual of a compact group is a discrete one (and vice versa).    
\snl
Other possible topics for further research in this field are the quantum double construction and other, similar methods to obtain new examples of algebraic quantum hypergroups. Also some nice applications would be most welcome.

\nl\nl

\bf References \rm
\nl
\parindent 1.0 cm
\item{[A]} E.\ Abe: \it Hopf algebras. \rm Cambridge University Press (1977).
\smallskip
\item{[C-V]} Y.A.\ Chapovsky and L.I.\ Vainerman: {\it Compact quantum hypergroups}. J. Operator Theory {\bf 41} (1999), 261-289.
\smallskip
\item{[DC-VD]} K.\ De Commer and A. Van Daele: {\it Multiplier Hopf algebras imbedded in C$^*$-algebraic quantum groups}. Preprint K.U.Leuven (2006).
\smallskip
\item{[D-VD1]} L.\ Delvaux and A.\ Van Daele: {\it Traces on (group-cograded) multiplier Hopf algebras}. Preprint  University of Hasselt and University of Leuven (2006).  
\smallskip
\item{[D-VD2]} L.\ Delvaux and A.\ Van Daele: {\it Algebraic quantum hypergroups. Examples and special cases}. Preprint  University of Hasselt and University of Leuven.  In preparation.
\smallskip
\item{[D-VD-W]} L.\ Delvaux, A.\ Van Daele and Shuanhong Wang: {\it A note on Radford's $S^4$ formula}. Preprint University of Hasselt, K.U.Leuven and Southeast University Nanjing (2006), math.RA/0608096.
\smallskip
\item{[Dr-VD-Z]} B.\ Drabant, A.\ Van Daele and Y.\ Zhang: {\it Actions of multiplier Hopf algebras}. Comm. Algebra {\bf 27} (1999), 4117-4127.
\smallskip
\item{[Ka]} A.A.\ Kalyuzhnyi: {\it Conditional expectations on quantum groups and new examples of quantum hypergroups}. Methods of Funct. Anal. Topol. {\bf 7} (2001), 49-68.
\smallskip
\item{[Ku]} J.\ Kustermans: {\it The analytic structure of algebraic quantum groups}. J.\ of Alg.\ {\it 259} (2003), 415--450.
\smallskip
\item{[K-V1]} J.\ Kustermans \& S.\ Vaes: \it A simple
definition  for locally compact quantum groups. \rm  C.R. Acad. Sci., Paris,
S\'er. I {\bf 328} (10) (1999), 871--876.
\smallskip
\item{[K-V2]} J.\ Kustermans \& S.\ Vaes: \it Locally compact
quantum groups. \rm Ann.\ Sci.\ \'Ec.\ Norm.\ Sup.\  {\bf 33} (2000), 837--934.
\smallskip
\item{[K-V3]} J.\ Kustermans \& S.\ Vaes: \it Locally compact
quantum groups in the von Neumann algebra setting. \rm Math.
Scand. {\bf 92} (2003), 68--92.
\smallskip
\item{[K-VD]} J.\ Kustermans \& A.\ Van Daele: {\it C$^*$-algebraic quantum groups arising from algebraic quantum groups}. Int.\ J.\ Math.\ {\bf 8} (1997), 1067--113.
\smallskip
\item{[L-VD1]} M.B.\ Lanstad and A.\ Van Daele: {\it Compact and discrete subgroups of algebraic quantum groups}.   Preprint University of Trondheim and University of Leuven (2006).
\smallskip
\item{[L-VD2]} M.B.\ Landstad and A.\ Van Daele: {\it Multiplier Hopf $^*$-algebras and groups with compact open subgroups}. Preprint University of Trondheim and University of Leuven (2006).
\smallskip
\item {[M-VD]} A.\ Maes \& A.\ Van Daele:  {\it Notes on compact
quantum groups}. Nieuw Archief voor Wiskunde, Vierde serie
{\bf 16} (1998), 73--112.
\smallskip
\item{[S]} M.\ Sweedler: {\it Hopf algebras}. Benjamin, New-York (1969). 
\smallskip
\item{[V]} L.I.\ Vainerman: {\it Gel'fand pair associated with the quantum groups of motions of the plane and $q$-Bessel functions}. Reports on Mathematical Physics {\bf 35} (1995), 303-326.
\smallskip
\item{[VD1]} A.\ Van Daele: {\it Multiplier Hopf algebras}. Trans. Am. Math. Soc. {\bf 342}(2) (1994), 917-932.
\smallskip
\item{[VD2]} A.\ Van Daele: {\it An algebraic framework for group duality}. Adv. in Math. {\bf 140} (1998), 323-366.
\smallskip
\item{[VD3]} A.\ Van Daele: {\it Multiplier Hopf $^\ast$-algebras with positive integrals: A laboratory for locally compact quantum groups.} Irma Lectures in Mathematical and Theoretical Physics 2: Locally compact Quantum Groups and Groupoids. Proceedings of the meeting in Strasbourg on Hopf algebras, quantum groups and their applications (2002). Ed.\ V.\ Turaev \& L.\ Vainerman. Walter de Gruyter, (2003), 229--247.
\smallskip
\item{[VD4]} A.\ Van Daele: {\it Locally compact quantum groups. A von Neumann algebra approach}. Preprint K.U.\ Leuven (2006), 49p. (math.OA/0602212).
\smallskip 
\item{[VD-W]} A.\ Van Daele and Shuanhong Wang: {\it The Larson Sweedler theorem for multiplier Hopf algebras}. J.\ of Alg.\ {\bf 296} (2006), 75--95.
\smallskip
\item{[VD-Z]} A.\ Van Daele and Y.\ Zhang: {\it A survey on multiplier Hopf algebras}. Proceedings of the conference in Brussels on Hopf Algebras and Quantum Groups.  Eds. Caenepeel/Van Oystaeyen (2000), 269-309.  Marcel Dekker (New York).
\smallskip
\item{[W1]} S.L.\ Woronowicz: {\it Compact matrix pseudogroups}. Comm. Math. Phys. {\bf 111} (1987), 613-665.
\smallskip
\item{[W2]} S.L.\ Woronowicz: {\it  Compact quantum groups.}
Quantum symmetries/Symm\'{e}tries quantiques.  Proceedings of the
Les Houches summer school 1995, North-Holland, Amsterdam (1998),
845--884.
\smallskip

\end

\item{[B-S]} S.\ Baaj \& G.\ Skandalis: \it Unitaires multiplicatifs
et dualit\'e pour les produits crois\'es de C$^*$-alg\`ebres. \rm
Ann. Scient.\ Ec.\ Norm.\ Sup., 4\`eme s\'erie,  {\bf 26} (1993), 425-488.
\smallskip

\item{[D-VD1]} L.\ Delvaux \& A.\ Van Daele: {\it The Drinfel'd double versus the Heisenberg double for an algebraic quantum group}. J.\ Pure \& Appl.\ Alg.\ 190 (2004), 59--84.
\smallskip

\item{[D-VD2]} L.\ Delvaux \& A.\ Van Daele: {\it Traces on (group-cograded) multiplier Hopf algebras}. Preprint K.U.Leuven and Universiteit Hasselt (2005).
\smallskip

\item{[E-S]} M.\ Enock \& J.-M.\ Schwartz: {\it Kac algebras and duality
for locally compact groups.} Springer (1992).
\smallskip

\item{[K]} E.\ Kirchberg: Lecture at the conference 'Invariance in Operator 
Algebras', Copenhagen, August 1992.
\smallskip

\item{[K-V4]} J.\ Kustermans \& S.\ Vaes: \it The operator algebra
approach to quantum groups. \rm Proc. Natl. Acad. Sci. USA {\bf 97 (2)}
(2000), 547--552.
\smallskip

\item{[K-*]} J.\ Kustermans, S.\ Vaes, L.\ Vainerman, A.\ Van Daele and S.L.\ Woronowicz: {\it Locally compact quantum
groups}. Lecture Notes School/Conference on Noncommutative Geometry and Quantum groups, Warsaw, 2001. Banach Centre
Publications, to appear.
\smallskip

\item{[M-VD2]} A.\ Maes \& A.\ Van Daele: {\it The multiplicative
unitary as a basis for duality.} Preprint K.U.\ Leuven (2002),
37p. (\#math.OA/0205284).
\smallskip

\item{[M-N]} M.\ Masuda \& Y.\ Nakagami: {\it A von Neumann algebra
framework for the duality of quantum groups.} Publ. RIMS Kyoto
{\bf 30} (1994), 799--850.
\smallskip

\item{[M-N-W]} T.\ Masuda, Y.\ Nakagami \& S.L.\ Woronowicz: {\it A C$^*$-algebraic framework for the quantum groups}. Int.\ J.\ of Math.\ {\bf 14} (2003), 903--1001.  
\smallskip

\item{[Q-VD]} J.\ Quaegebeur \& A.\ Van Daele: {\it Weights on C$^*$-algebras}. Preprint K.U.Leuven. In preparation.
\smallskip

\item{[R]} D.E.\ Radford: {\it The order of the antipode of a finite-dimensional Hopf algebra is finite}. Amer.\ J.\ Math.\ {\bf 98} (1976), 333--355.
\smallskip
 
\item{[S-W]} P.\ Soltan \& S.L.\ Woronowicz: {\it A remark on manageable multiplicative unitaries}. Letters on Math.\ Phys. {\bf 57} (2001), 239--252.
\smallskip

\item{[St]} S.\ Stratila: {\it Modula Theory in Operator Algebras}. Abacus Press, Turnbridge Wells, 1981.
\smallskip

\item{[S-V-Z1]} S.\ Stratila, D.\ Voiculescu \& L.\ Zsido: {\it Sur les produits crois\'es}. Comptes Rendus Acad.\ Sc.\ Paris {\bf 282} (1976), ??--??
\smallskip

\item{[S-V-Z2]} S.\ Stratila, D.\ Voiculescu \& L.\ Zsido: {\it On crossed products I}. Rev.\ Roumaine Math.\ Pure et Appl. {\bf 21} (1976), 1411--1449.   
\smallskip

\item{[S-V-Z3]} S.\ Stratila, D.\ Voiculescu \& L.\ Zsido: {\it On crossed products II}. Rev.\ Roumaine Math.\ Pure et Appl. {\bf 22} (1977), 83--117.   
\smallskip

\item{[S-Z]} S.\ Stratila \& L.\ Zsido: {\it Lectures on von Neumann algebras}. Abacus Press, Tunbridge
Wells, 1979.
\smallskip

\item{[T1]} M.\ Takesaki: {\it Tomita's theory of modular Hilbert algebras and applications}. Springer Lecture Notes in Mathematics 128. Springer (1970).
\smallskip

\item{[T2]} M.\ Takesaki: {\it Theory of Operator Algebras I}. Springer-Verlag, New York (1979). 
\smallskip

\item{[T3]} M.\ Takesaki: {\it Theory of Operator Algebras II}. Springer-Verlag, New York (2001).
\smallskip

\item{[Tay]} D.C.\ Taylor: {\it The strict topology for double centralizer algebras}. Trans.\ Am.\ Math.\ Soc.\ {\bf 150} (1970), 633--643.
\smallskip

\item{[V1]} S.\ Vaes: {\it A Radon-Nikodym theorem for von Neumann algebras}. J.\ Operator Theory {\bf 46} (2001), 477-489.
\smallskip

\item{[V2]} S.\ Vaes: {\it Locally compact quantum groups}. PhD. Thesis K.U.Leuven (2001).
\smallskip

\item{[V-VD]} S.\ Vaes \& A.\ Van Daele: {\it Hopf C$^*$-algebras}.
Proc.\ London Math.\ Soc.\ {\bf 82} (2001), 337--384.
\smallskip

\item{[VD1]} A.\ Van Daele: {\it Dual pairs of Hopf $^*$-algebras}. Lecture notes of a talk at the University of 
Orl\'eans. Preprint K.U.Leuven (1991).
\smallskip
 
\item{[VD2]} A.\ Van Daele: {\it Multiplier Hopf algebras.} Trans.\
Amer.\ Math.\ Soc.\ {\bf 342} (1994), 917--932.
\smallskip

\item{[VD3]} A.\ Van Daele: {\it An algebraic framework for group duality.}
Adv.\ in Math.\ {\bf 140} (1998), 323--366.
\smallskip

\item{[VD4]} A.\ Van Daele: {\it Quantum groups with invariant integrals}. Proc.\ Natl.\ Acad.\ Sci.\ USA {\bf 97} (2000), 541-546.
\smallskip

\item{[VD5]} A.\ Van Daele: {\it The Haar measure on some locally compact quantum groups}. Preprint K.U.\ Leuven (2001)., math.OA/0109004.
\smallskip

\item{[VD6]} A.\ Van Daele: {\it Locally compact quantum groups: The von Neumann algebra versus the C$^*$-algebra approach}. Preprint K.U.\ Leuven (17p) (2005). To appear in The Bulletin of Kerala Mathematics Association.
\smallskip

\item{[VD7]} A.\ Van Daele: {\it Notes on Multiplier Hopf ($^*$-)algebras and algebraic quantum groups}. Preprint K.U.\ Leuven. In preparation.
\smallskip
 
\item{[VD8]} A.\ Van Daele: {\it Notes on Locally Compact Quantum Groups}. Preprint K.U.\ Leuven. In preparation.
\smallskip
 
\item{[VH]} L.\ Vanheeswijck: {\it Duality in the theory of
crossed products}. Math.\ Scand.\ {\bf 44} (1979), 313--329.
\smallskip

\item{[W2]} S.L.\ Woronowicz: {\it From multiplicative unitaries to quantum groups}. Int.\ J.\ Math.\ {\bf 7} (1996), 127--149.

\end